\documentclass[notitlepage]{scrartcl}
\usepackage[shortlabels]{enumitem}
\usepackage{amsmath}
\usepackage{amssymb}
\usepackage{amsthm}
\usepackage{abstract}
\usepackage{xcolor}
\usepackage{comment}
\usepackage{mathtools}
\usepackage{graphicx}
\usepackage{caption}
\usepackage{subcaption}
\usepackage{float}
\usepackage{hyperref}
\usepackage{enumitem}
\makeatletter
\def\namedlabel#1#2{\begingroup
   \def\@currentlabel{#2}%
   \label{#1}\endgroup
}
\makeatother
\usepackage{dsfont}
\hypersetup{
    colorlinks=true,
    linkcolor=blue,
    filecolor=magenta,      
    urlcolor=cyan
    }
\makeatletter
\newcommand*{\barfix}[2][.175ex]{%
  \mathpalette{\@barfix{#1}}{#2}%
}
\newcommand*{\@barfix}[3]{%
  \vbox{%
    \kern#1\relax
    \hbox{$#2#3\m@th$}%
  }%
}
\makeatother

\usepackage[capitalise]{cleveref}
\newtheorem{theorem}{Theorem}
\newtheorem{thm}{Theorem}[section]
\newtheorem{corollary}[thm]{Corollary}
\newtheorem{lemma}[thm]{Lemma}
\newtheorem{proposition}[thm]{Proposition}
\newtheorem{claim}[thm]{Claim}

\newcommand{\brac}[1]{\left( #1 \right)}
\newcommand\bfrac[2]{\left(\frac{#1}{#2}\right)}

\theoremstyle{definition}

\newtheorem{definition}[thm]{Definition}
\newtheorem{remark}[thm]{Remark}

\newcommand{\footremember}[2]{%
    \footnote{#2}
    \newcounter{#1}
    \setcounter{#1}{\value{footnote}}%
}

\usepackage[margin=1in]{geometry} 

\def\cB{\mathcal{B}}

\def\cP{\mathcal{P}}

\def\cE{\mathcal{E}}

\def\cM{\mathcal{M}}

\def\eps{\varepsilon}

\title{\vspace{-1.5cm} Linear-length cycles in the supercritical percolated hypercube} 
\author{%
Michael Anastos \footremember{alley1}{\scriptsize{Department of Mathematics, ETH, Zürich, Switzerland. Email: michael.anastos@math.ethz.ch}}}
\date{\vspace{-1.5cm}}

\begin{document}
\maketitle

\begin{abstract}
Let $Q_p^d$ be the random subgraph of the $d$-dimensional binary hypercube obtained by retaining each edge independently with probability $p$. A long-standing folklore conjecture states that the emergence of the giant component in $Q_p^d$, at $p\approx 1/d$, is accompanied by the emergence of cycles of linear length. We confirm this conjecture and determine the order of the longest cycle near the phase transition. Specifically, for every sufficiently small constant $\varepsilon>0$, we show that if $p=(1+\varepsilon)/d$, then with high probability the longest cycle of $Q_p^d$ has length $\Theta(\varepsilon^2 2^d)$. The dependence on $\varepsilon$ is optimal up to constant factors.
 
\end{abstract}

\section{Introduction}

The \emph{$d$-dimensional binary hypercube} $Q^d$ is the graph with vertex set $\{0,1\}^d$, in which two vertices are adjacent if and only if they differ in a single coordinate. The \emph{percolated hypercube} $Q^d_p$ is the random subgraph of $Q^d$ obtained by retaining each edge independently with probability $p$. The hypercube and its random subgraphs have been studied extensively in combinatorics, probability and computer science, where they arise in a wide range of problems.

The study of the percolated hypercube $Q^d_p$ was pioneered by Burtin~\cite{B77} and Sapo\v{z}enko~\cite{S67}, who showed that $p=1/2$ is a sharp threshold for connectivity.  Results on its component structure followed. Erd\H{o}s and Spencer~\cite{ES79}  observed that, if $p=c/d$ with $c<1$, then $Q^d_p$ typically contains only components of order $O(d)$. On the other hand, Ajtai, Koml\'os and Szemer\'edi~\cite{AKS81} proved that, if $p=c/d$ with $c>1$, then $Q^d_p$ typically contains a connected component $L_1$ with linearly many vertices, called the \emph{giant component}. These results established that the size of the largest component of $Q_p^d$ goes through a phase transition at $p=(1+o(1))/d$.

Bollob\'as, Kohayakawa and \L{}uczak~\cite{BKL92} determined the size of the giant more precisely. For $c>0$, let $y(c)$ denote the survival probability of a Galton--Watson process with offspring distribution $\operatorname{Poisson}(c)$. They proved that, if $p=c/d$ with $c>1$, then typically $L_1$ has $y(c)2^d+o(2^d)$ vertices, while every other component has order $O(d)$. Thus, in the corresponding regime, the component structure of $Q^d_p$ resembles that of the binomial random graph $G(n,p)$~\cite{ER60}, exhibiting the so-called Erd\H{o}s--R\'enyi component phenomenon~\cite{DK23,DK25}.

A basic question in the study of random graphs is whether the emergence of the giant component is accompanied by the emergence of cycles of linear length.
In this paper we address this question for the percolated hypercube. The circumference of a graph is the length of its longest cycle. We prove that for every sufficiently small $\eps>0$ and $p=(1+\eps)/d$, the circumference of $Q_p^d$ is of order $\eps^22^d$ with high probability.

\subsection{The circumference of $G(n,p)$}
Erd\H{o}s~\cite{erdos74} conjectured that for every constant $c>1$, $G(n,c/n)$ has a cycle of linear length whp. His conjecture was resolved by Ajtai, Koml\'os and Szemer\'edi~\cite{AKS81a}, who showed that for any constant $c>1$, if $p\ge c/n$, then  $G(n,p)$  contains a cycle of linear length whp. A weaker version of this result, covering the range $c>2\log 2$, was proved earlier by Fernandez de la Vega~\cite{Fdlv79}.  These results also give that the fraction of vertices missed by a longest cycle typically tends to zero as $c$ tends to infinity. Later, Ajtai, Komlós and Szemerédi~\cite {AKS85}, and independently Bollobás~\cite {B84}, determined that the threshold for the existence of a Hamilton cycle in the binomial random graph equals the threshold for having minimum degree at least $2$, completing the picture.

\subsection{The circumference of $Q_p^d$}

A folklore conjecture states that the circumference of $Q^d_p$ exhibits behaviour analogous to that of $G(n,p)$. Namely, it undergoes a phase transition near $p=1/d$, from logarithmic to linear in the number of vertices; it then continues to grow with the giant component, and finally $Q^d_p$ becomes Hamiltonian as soon as its minimum degree reaches $2$.

A breakthrough came from Condon, Espuny D\'iaz, Gir\~ao, K\"uhn and Osthus~\cite{CDGKO21}, who proved
the Hamiltonicity part of the conjecture. In particular, they showed that $p=1/2$ is the threshold for the existence of a Hamilton cycle in $Q^d_p$.

A second breakthrough came from the author together with Diskin, Erde, Kang, Krivelevich and Lichev~\cite{anastos2025nearly}, who proved the ``middle'' part of the conjecture in a strong form. Namely, they showed that, for every $\eps\in(0,1)$, there exists a constant $c=c(\eps)$ such that, if $p\ge c/d$, then $Q^d_p$ typically contains a cycle covering at least $(1-\eps)2^d$ vertices.

In this paper, we establish the emergence of cycles of linear length
just above $p=1/d$, showing that the circumference of $Q^d_p$
indeed undergoes a phase transition near this point
(see Theorem~\ref{thm:main}). This question has been reiterated over the years; see, for example, the 2014 ICM survey by Frieze~\cite{F14}, \cite[Question 5.5]{EKK22}, \cite[Question 72]{F19}, \cite[Conjecture 1.2]{anastos2025nearly}.
In this regime, by exploiting the expansion properties of linear-sized
sets, Diskin, Erde, Kang and Krivelevich~\cite{DEKK24} proved that, if $p=c/d$ for a constant $c>1$, then $Q^d_p$ typically contains a cycle of
length $\Omega(2^d/(d\log d))$. We prove that the cycle can be taken
to have linear length.
\begin{theorem}\label{thm:main}
Let $\varepsilon>0$ be a sufficiently small constant. There exists
a constant $C>0$ that is independent of $\eps$ such that, if $p=p(d)=(1+\varepsilon)/d$, then with high probability\footnote{We say that a sequence of events $(\cE_d)_{d\geq 1}$ holds \emph{with high probability} (whp for short) if $\Pr(\cE_d)\to 1$ as $d\to \infty$. } $Q^d_p$
contains a cycle of length at least $C \eps^2 2^d$.
\end{theorem}
The next theorem shows that the dependence on $\eps$ is optimal. The $2$-core of a graph $G$ is the maximal subgraph of $G$ of minimum degree at least $2$. Every cycle of $G$ lies in its $2$-core. Thus the circumference of $G$ is at most the size of the largest component of its $2$-core.  
\begin{theorem}\label{thm: size of 2-core}
Let $\eps>0$ be a sufficiently small constant and let $p=p(d)=(1+\eps)/d$. Then, whp the $2$-core of the giant component of $Q_p^d$ has size $(2\eps^2+O(\eps^3))2^{d}$.
\end{theorem}

Theorem 1 settles the existence question. Theorem 2 shows why the scale $\eps^22^d$ is the right one near the transition. Together, the two theorems show that, for $p=(1+\eps)/d$, where $\eps>0$ is a sufficiently small constant, the circumference of $Q_p^d$ is $\Theta(\eps^2 2^d)$ whp.

Our proof uses the expansion properties of medium-sized sets and their resilience under vertex deletion. In addition, we initiate the study of structural properties of the $2$-core of $Q_p^d$. We believe that these ideas could also be useful in future work.

\subsection{Proof outline}\label{sec:outline}

Our proof combines three ingredients: a path-packing result for locally expanding graphs, a method for joining the resulting paths through successive subcubes, and expansion estimates for medium-sized sets in the $2$-core.

We choose pairwise vertex-disjoint subcubes $G'_0$ and $G'_1,\ldots,G'_k$ of $Q^d$. The initial subcube has dimension $d-2$, while each of the remaining subcubes has dimension $d-O_\varepsilon(1)$. Consecutive subcubes $G'_{j-1}$ and $G'_j$ are joined by a perfect matching $M_j$ in $Q^d$ for $j\geq 2$, and $G_0'$ is joined to $G_1'$ by a set of pairwise disjoint paths. For $0\le j\le k$, let $G_j$ be the subgraph of $Q_p^d$ induced by $V(G'_j)$.

We first find vertex-disjoint paths of at least medium length in the $2$-core of the giant component of $G_0$ that cover a linear number of vertices. We then join these paths into a single cycle using the subsequent graphs $G_1,\ldots,G_k$ and the matching edges retained in $Q_p^d$. Here, ``medium'' means exponential in $d$, with a small positive exponent depending on $\eps$. Working in the $2$-core at the packing stage gives the optimal dependence on $\eps$.

\medskip
\noindent\textbf{Finding the paths.}
This part of the proof uses a variant of a result of Krivelevich~\cite{krivelevich2019long}, which gives a long cycle from local expansion. We take a maximal family of vertex-disjoint cycles of at least medium length and remove their vertices. The remaining graph contains no further cycle of this length. The same result then implies that any component larger than the chosen scale contains a poorly expanding almost connected set of medium size. We successively remove such sets until all remaining components are small. Removing one edge from each of the original cycles gives the required paths.

The goal is now to show that the vertices outside the paths cannot greatly outnumber those on the paths. Each set removed in the second stage expands in the original graph, but has only a small neighbourhood in the graph remaining at the time of its removal. Thus, many of its neighbours must lie on the paths or in sets removed earlier. Using the accompanying bounds on the number of edges incident with these sets, a double-counting argument shows that the sets removed in this stage together contain at most a constant multiple of the number of vertices on the paths.

It remains to control the components left after these removals. We divide them into those with good expansion and those with poor expansion, measured in the original graph. All their external neighbours lie on the paths or in the removed sets. For the components with good expansion, counting the edges between these components and the removed vertices bounds their total order in terms of the number of removed vertices. On the other hand, too many poorly expanding components would allow us to combine a suitable collection of them with a segment of a single cycle or a single set of removed vertices to form an almost connected set. Most of its vertices would come from the poorly expanding components. The upper neighbourhood bound then gives a poorly expanding set, contradicting expansion in the original graph. Together, these comparisons show that the paths cover a positive proportion of the base graph. Applied to the $2$-core, this gives the paths needed for the next step.

\medskip
\noindent\textbf{Joining the paths.}
Let $P_1,\ldots,P_t$ be the resulting paths. For each $P_i$, we choose sets $S_i^+$ and $S_i^-$ in its first and last thirds, respectively. We choose these sets so that each of their vertices is connected to its projection in the first auxiliary subcube along a retained path, and then to its projection in the subsequent subcubes along retained matching edges. In addition, its projection in each $G_j$ belongs to the giant component of $G_j$. We may take all the sets to have the same size $s$, which is exponentially large in $d$. Their projections into the subsequent subcubes form the terminal sets through which we join the original paths.

In each round, we work in the next giant component and choose a maximal family of vertex-disjoint short paths joining the available terminal sets. These paths avoid the unused terminal sets and do not close a cycle before all the original paths have been joined. After deleting their vertices, we divide the terminal sets still to be joined into two roughly equal groups, and let $U_1$ and $U_2$ be their unions. If sets of size at least $\min\{|U_1|,|U_2|\}$ still expand sufficiently in the remaining graph, then we can find another short path between $U_1$ and $U_2$, contradicting maximality.

Thus the number of vertices in terminal sets still to be joined is controlled by the number of vertices already used. Comparing this deletion cost with the expansion of the giant component shows that, if all relevant sets of size at least $\min\{|U_1|,|U_2|\}$ initially expand by a factor of, say, $1/d^2$, then the number of deleted vertices is larger than $\min\{|U_1|,|U_2|\}/2d^2$; otherwise, enough expansion remains to find another connection. Consequently, at each round the number of remaining terminal sets decreases  by a factor of order $s/d^{O(1)}$. Since $s$ is exponentially large in $d$, only $O(d/\log s)$ rounds are needed. Choosing $k$ sufficiently large, we obtain a cycle retaining a third of the vertices on the original paths, and hence of linear length.

\medskip
\noindent\textbf{Why the $2$-core.} 
The fraction of vertices covered in the first step, and hence retained on the cycle, is bounded below by a multiple of the ratio between the lower bound on vertex expansion and the upper bound on edge expansion of medium-sized almost connected sets. Taking the base graph to be the giant of $Q_p^d$ gives a cycle covering a $\Theta(\varepsilon^2)$ fraction of the vertices of $L_1$, and thus of length $\Theta(\varepsilon^3 2^d)$ (see \cite{DEKK24,anastos2025diameter} for the expansion properties of $L_1$).

To obtain the optimal dependence on $\varepsilon$, we instead work with the $2$-core of the giant. Theorem~\ref{thm: expansion of 2-core} shows that the corresponding ratio is bounded below by a constant. Thus, the resulting cycle retains a constant proportion of the vertices of the $2$-core. Then, Theorem~\ref{thm: size of 2-core} gives that the $2$-core typically has size $\Theta(\varepsilon^2 2^d)$, yielding a cycle of the desired length $\Theta(\varepsilon^2 2^d)$.

\subsection{Plan of the paper}
In Section~2, we introduce the notation and collect the preliminary results. Section~3 states the main ingredients and uses them to prove Theorem~\ref{thm:main}. Sections~4 and~5 prove the packing and joining results, respectively. Section~6 establishes the size and expansion estimates for the $2$-core, including Theorem~\ref{thm: size of 2-core}.

\section{Notation and Preliminaries}\label{sec:prel}

Given a graph $G=(V,E)$, for $A,B\subseteq V$, we write $e_G(A,B)$ for the number of edges with one endpoint in $A$ and one in $B$, with edges inside $A\cap B$ counted twice. We write $N_G(A)$ for the set of vertices in $V\setminus A$ with a neighbour in $A$, $G[A]$ for the subgraph induced by $A$, and $e_G(A):=|E(G[A])|$. For $S\subseteq V(G)$, we say that $S$ is connected in $G$ if $G[S]$ is connected. We say that $S$ is \emph{almost connected in $G$} if it is connected in $G$, or if $G[S\cup\{v\}]$ is connected for some $v\in V(G)\setminus S$. A path $P\subseteq G$ is \emph{bare} in $G$ if all of the internal vertices of $P$ have degree two in $G$.

We omit the dependence on the base graph when it is clear from the context.

We also use the giant-component result recalled in the introduction: for fixed $c>1$, the largest component of $Q_{c/d}^d$ has $(y(c)+o(1))2^d$ vertices, and every other component has order $O(d)$, whp~\cite{BKL92}. Here $y(c)>0$ is the survival probability of a Galton--Watson process with offspring distribution $\operatorname{Poisson}(c)$.  
 
Throughout the paper, we omit rounding signs whenever this does not affect the validity of our arguments. We use standard asymptotic notation, with subscripts for dependencies of implicit constants. All logarithms without a subscript are natural logarithms. We also assume that $d$ and $1/\eps$ are sufficiently large for the calculations below to hold.

\subsection{Auxiliary graph-theoretic results}
We start with a bound on the number of rooted trees in a graph of given maximum degree.

\begin{lemma}[\cite{beveridge1998random}, Lemma 2]\label{lem: tree count} Let $G$ be a graph with maximum degree $\Delta$, let $v\in V(G)$ and $k\geq 1$ be an integer. Let $t(v,k)$ be the number of rooted trees in $G$ which have root $v$ and $k$ vertices. Then,
\[
 t(v,k) \leq \frac{k^{k-2}\Delta^{k-1}}{(k-1)!}\leq (e\Delta)^{k-1}.\]
\end{lemma}

Next, we state Harper's edge-isoperimetric inequality for hypercubes,  see~\cite{H64} and also \cite{B67,H76,L64}.
\begin{lemma}\label{lem: Harper}
For every integer $d \ge 1$ and for every set $S\subseteq V(Q^d)$,
\begin{align*}
e(S,V(Q^d)\setminus S)\ge |S|\left(d-\log_2|S|\right).
\end{align*}
\end{lemma}

\subsection{Concentration inequalities}
We use the following version of the Chernoff bound for binomial random variables (see, for example, Section 21.4 of \cite{frieze1996introduction}).
\begin{theorem}[Chernoff bound]\label{thm: chernoof}
For a binomial random variable $X$ and any $t \in [0,\mathbb{E}[X]]$,
\[
\Pr[|X-\mathbb E X|\ge t]\le 2\exp\left(-\frac{t^{2}}{3 \mathbb{E}[X]}\right).
\]
\end{theorem}
We will also use McDiarmid's  inequality, stated below (see Theorem 1 of \cite{McDiarmid:89}). 
 \begin{theorem}\label{thm:mcdiarmid}
Let $n\geq1$ be an integer, $C>0$, and $f:\{0,1\}^n \to \mathbb{R}$ be such that if $x,x'\in \{0,1\}^n$ differ only in a single coordinate, then $|f(x)-f(x')|\le C$.
Let $X_{1},\dots, X_{n}$ be independent Bernoulli random variables. Let  $X=f(X_1,\dots,X_n)$. Then, for any $t\geq 0$ we have that
\[
\Pr[|X-\mathbb E X|\ge t]\le 2\exp\left(-\frac{2t^{2}}{C^2n}\right).
\]
\end{theorem}

\subsection{Auxiliary results on percolation}
We first recall the standard approximation to $y(1+\eps)$, the survival probability of a Galton--Watson process with offspring distribution $\operatorname{Poisson}(1+\eps)$, for sufficiently small $\eps>0$ (see, for example, \cite{haldane1927mathematical}). The same approximation applies to offspring distribution $\operatorname{Binomial}(d,(1+\eps)/d)$ for sufficiently large $d$.
\begin{lemma}\label{lem: approximation survival probability}
Let $\eps>0$ be a sufficiently small constant. Then, $y(1+\eps)=2\eps+O(\eps^2)$.  
\end{lemma}
The following lemma records expansion properties of the giant component.
\begin{lemma}\label{lem: giant weak expansion}
Let $\eps>0$ be a sufficiently small constant and $p=p(d)=(1+\eps)/d$. Whp the giant component $L_1$ of $Q^d_p$ has order  $(2\eps+O(\eps^2))2^{d}$ and the following property. For every set $S\subseteq V(L_1)$ of size $|S|\leq 2|V(L_1)|/3$ we have $|N_{L_1}(S)|\geq |S|/d^2$. Finally, whp the second largest component of $Q_p^d$ has order less than $d^2$. 
\end{lemma}
 \begin{proof}
The bound on the second largest component and the estimate on $|V(L_1)|$ follow from \cite{BKL92}. For the expansion bound, \cite[Corollary 6.3]{anastos2025diameter} implies that whp for every set $S\subseteq V(L_1)$ of size $|S|\leq 2|V(L_1)|/3$ we have $e_{L_1}(S,V(L_1)\setminus S)\geq |S|/(d\log\log d)$. The probability that $Q_p^d$ has a vertex of degree larger than $t=2d/\log d$ is at most  $2^d \binom{d}{t} p^{t} \leq 2^d (epd/t)^{t} = 2^d e^{-(1+o(1))t\log t}=o(1)$. Combining these two bounds gives that whp for every set $S\subseteq V(L_1)$ of size $|S|\leq 2|V(L_1)|/3$ we have $|N_{L_1}(S)|\geq (|S|/(d\log\log d))/(2d/\log d)\geq |S|/d^2$, which completes the proof.  
 \end{proof}
The expansion bound in the lemma above also controls the components outside the $2$-core and the bare paths within it.
\begin{corollary}\label{cor: sizes of paths bare paths}
Let $\eps>0$ be a sufficiently small constant and $p=p(d)=(1+\eps)/d$. Let $C_2$ be the $2$-core of the giant component $L_1$ of $Q_p^d$. Then, whp every component of $Q_p^d-V(C_2)$ has at most $d^3$ vertices and every bare path of $C_2$ has at most $d^3-3$ vertices.  
\end{corollary}
\begin{proof}
First, suppose that $Q_p^d-V(C_2)$ contains a component $C$ of order larger than $d^3$. Then, Lemma~\ref{lem: giant weak expansion} implies that $C\subseteq L_1$. Since $C$ does not intersect $C_2$, it is acyclic and $N_{L_1}(C)\subseteq C_2$ contains at most one vertex. Therefore, there exists a subtree $T$ of $C$ of order $|V(T)|\in [d^3,d^4+1]$ that has at most one neighbour in $L_1$. Lemma~\ref{lem: giant weak expansion} implies that this occurs with probability $o(1)$.

Now suppose that $C_2$ contains a bare path on at least $d^3-3$ vertices. Let $C$ be the subgraph of $L_1$ consisting of $d^3-5$ consecutive internal vertices along with the trees in $L_1-V(C_2)$ attached to them. In the event that no such tree has order larger than $d^3$,  $C$ has order in $[d^3-5, d^7]$ and has at most two external neighbours in $L_1$ (which also belong to $C_2$). Lemma~\ref{lem: giant weak expansion} implies that this occurs with probability $o(1)$.    
\end{proof}
The next lemma estimates the probability that a given vertex belongs to a sufficiently large component, and that its neighbourhood reaches a prescribed size.
\begin{lemma}
Let $\eps>0$ be sufficiently small and $p=p(d)=(1+\eps)/d$. Let $r\in V(Q^d)$. Then,
\begin{equation}\label{eq: exploring to medium}
       \Pr(\text{the component of  $r$ in $Q_p^d$ has at least $\eps^{-6}$ vertices}) = 2\eps+O(\eps^2).
   \end{equation}
and
\begin{equation}\label{eq: expansion}
       \Pr\brac{\text{there are at least $d^{20}$ vertices at distance at most $\log^2 d$ from $r$ in $Q_p^d$}} = 2\eps+O(\eps^2),
\end{equation}
\end{lemma}

\begin{proof} 
Let $T_r$ be the component of $r$ in $Q_p^d$. Let $\eps^{-6}\leq i< d^2$. By Lemma~\ref{lem: tree count}, the number of $i$-vertex trees in $Q^d$ containing $r$ is at most $(ed)^{i-1}$. Lemma~\ref{lem: Harper} gives that a connected set on $i$ vertices spans at most $(i \log_2 i)/2$ edges, so it has at least $(1+o(1))di$ edges leaving it. Hence, 
\begin{align*}
    \Pr(|V(T_r)|=i)&\leq (ed)^{i-1}p^{i-1}(1-p)^{(1+o(1))di} \leq \brac{epde^{-(1+o(1))pd}}^i 
    \\&\leq \exp \brac{\brac{1+\log (1+\eps)-(1+o(1))(1+\eps)}i}
    \leq  \exp \brac{-\eps^2i/10}.
\end{align*}
It follows that $|V(T_r)|\in [\epsilon^{-6}, d^2)$ with probability at most $O(\epsilon^2)$. Lemma~\ref{lem: giant weak expansion} implies that $T_r$ has at least $d^2$ vertices with probability $2\eps+ O(\epsilon^2)$. \eqref{eq: exploring to medium} follows. 

To prove \eqref{eq:expansion}, for $i\geq 0$ let $S_i$ be the set of vertices at distance $i$ from $r$ in $Q_p^d$, and let $S_{<i}\coloneqq \bigcup_{j=0}^{i-1}S_j$. In the event that $|S_i|\leq \eps^4|S_{<i}|$, and $\log^{3/2}d\leq |S_{<i}|\leq d^{20}$, there is a tree rooted at $r$ (of depth $i$) in $Q_p^d$ with $\log^{3/2}d\leq k\leq 2d^{20}$ vertices, of which at most $2\eps^4 k$ vertices have a neighbor outside this set. Lemma~\ref{lem: Harper} implies that its vertex set spans at most $(k\log_2k)/2$ edges of $Q^d$. By Lemma~\ref{lem: tree count}, this occurs with probability at most
\begin{align*}
&\sum_{k=\log^{3/2}d }^{2d^{20}}  (ed)^{k-1}p^{k-1}\binom{k}{2\eps^4 k}(1-p)^{(k-2\eps^4 k)d-k\log_2 k}
\\& \leq \sum_{k\geq \log^{3/2}d } \brac{ (edp)\bfrac{e}{2\eps^4 }^{2\eps^4}e^{-(1-O(\eps^4))dp}}^k  \leq \sum_{k\geq \log^{3/2}d } \brac{\bfrac{e}{2\eps^4}^{2\eps^4}e^{1+\log (1+\eps)-(1-O(\eps^4))(1+\eps)}}^k
\\& \leq \sum_{k\geq \log^{3/2}d } \brac{\bfrac{e}{2\eps^4}^{2\eps^4}e^{-\eps^2/4}}^k = \sum_{k\geq \log^{3/2}d } \brac{\bfrac{e}{2\eps^4}^{2\eps^2}e^{-1/4}}^{\eps^2 k}=o(d^{-1}).
\end{align*}
Lemma~\ref{lem: giant weak expansion} implies that $|S_{<\log^{3/2}d}|\geq \log^{3/2}d$ with probability at least $2\eps+O(\eps^2)$. The preceding bound gives 
\[|S_{<\log^2 d}|\geq \min\{d^{20},(1+\eps^4)^{\log^2d-\log^{3/2}d}(\log^{3/2}d)\}= d^{20}\] with probability at least $2\eps+O(\eps^2)$. On the other hand Lemma~\ref{lem: giant weak expansion} implies that $T_r$ has at most $d^2$ vertices with probability at least $1-2\eps+O(\eps^2)$, hence it has at least $d^{20}$ vertices with probability at most $2\eps+O(\eps^2)$. This completes the proof of \eqref{eq: expansion}.
\end{proof}

The last lemma of the section is part of \cite[Lemma 7.7]{anastos2025diameter} and will be used later in the proof of \cref{thm:main}.
\begin{lemma}\label{lem: giant is well spread} 
Let $c>1$ be constant and let $p=c/d$. Whp every vertex $v\in V(Q^d)$ is within distance two in $Q^d$ from a vertex of the giant component $L_1$ of $Q_p^d$.
\end{lemma}

\section{Proof of Theorem~\ref{thm:main}}
We use three results, whose proofs are given in Sections~\ref{sec: path packing}-~\ref{sec: proof expanstion}. The first yields a family of vertex-disjoint cycles, the second joins the resulting paths into a cycle, and the third gives the expansion estimates for the $2$-core needed to obtain the optimal dependence on $\eps$ and $d$.

Theorem~\ref{thm: cycle packing} implies, by removing one edge from each cycle, that if every medium-sized almost connected set in a graph $G$ has suitable expansion properties, then there exists a family of pairwise vertex-disjoint medium-length paths that covers a large proportion of $V(G)$.
\begin{theorem}\label{thm: cycle packing}
Let $k\geq 1$ be an integer, $0<c\leq C$, and $I$ be the interval $[I_{\min}, I_{\max}]$, where $I_{\min}=\min\{ ck/10,  ck/10C\}$ and $I_{\max}= \max\{10Ck/c,10Ck\}$. Let $G$ be a connected graph on more than $k$ vertices such that every almost connected set $W\subseteq V(G)$ with $|W|\in I$ satisfies the following properties:

\vbox{
\begin{enumerate}[label = \upshape{(P\arabic*)}, leftmargin= \widthof{R100000}]
    \item $c|W|\leq |N(W )|\leq C|W|$, and\label{item: vertex expansion}
    \item for every $Z\subseteq N(W)$ we have that  $e(W,Z)\leq |Z|+c^2k/10$.  \label{item: vertex implies edge}
\end{enumerate}
}
Let $\ell\coloneqq  ck/10 $. If $2 < 
\ell $, then $G$ contains a family of pairwise vertex disjoint cycles of length at least $\lceil \ell \rceil$, covering at least 
\[\frac{1}{7}\left(\frac{c}{C}\right) |V(G)|\]
vertices.
\end{theorem}

In the proof of Theorem~\ref{thm:main}, we work with a sequence of vertex-disjoint subcubes $H_0,\dots H_k\subseteq Q^d$. We first apply Theorem~\ref{thm: cycle packing} to find vertex-disjoint paths in the percolated subgraph of $H_0$. We then connect them successively through $H_1,\dots H_k$ until we obtain a single cycle covering at least a third of their vertices. The connections use projections of selected vertices in the first and last thirds of the paths. We call these projected sets \emph{terminal sets}. Proposition~\ref{lem:connecting} implies that this is possible when the graph sequence is sufficiently long and the graphs used for joining have suitable expansion properties. To keep track of the length retained as paths are joined, we use marked trees. Given a graph $G$, a marked tree is a triple $B=(T,S_1,S_2)$ where $T\subseteq G$ is a tree, $S_1,S_2\subseteq V(T)$ are non-empty, and $S_1\cap S_2=\emptyset$. We call $S_1,S_2$ the terminal sets of $B$. The length of a marked tree $B=(T,S_1,S_2)$, denoted by $d(B)$, is defined to be 
\[
  d(B):=\min\{d_T(x,y):x\in S_1,\ y\in S_2\}.
\]
We say that the marked trees $B=(T,S_1,S_2)$ and $B'=(T',S_1',S_2')$ are disjoint if the trees $T$ and $T'$ are vertex disjoint.

\begin{proposition}\label{lem:connecting}
Let $n,\ell,k,r,d,s\geq 1$ be integers. Let $G$ be an $n$-vertex graph, and let $G_0,G_1,\dots,G_k$ be pairwise vertex-disjoint subgraphs of $G$.  For $*\in\{+,-\}$, $i\in[r]$ and $j\in\{0,1,\dots,k\}$, let $S_{i,j}^*\subseteq V(G_j)$ be of size $|S_{i,j}^*|=s$, and assume that all these sets are pairwise disjoint.  Suppose that the following conditions hold.

\begin{enumerate}[label=\textnormal{(\alph*)},leftmargin=2.7em]
\item \label{item: distance trees} There are pairwise vertex-disjoint trees
      $T_1,\dots,T_r\subseteq G_0$ such that for $i\in [r]$ we have that $B_i^0=(T_i,S_{i,0}^+,S_{i,0}^-)$ is a marked tree of length at least $\ell$.
\item \label{item:  perfect matchings between levels}  $G$ contains a perfect matching $M_{i,j}^*$ between $S_{i,j}^*$ and $S_{i,j+1}^*$, for every $*\in\{+,-\}$,  $i\in[r]$, and $j\in\{0,\dots,k-1\}$.
\item \label{item: expansion2} For every $j\in[k]$ and every $A\subseteq V(G_j)$ if $|A|\le 2|V(G_j)|/3$, then $|N_{G_j}(A)|\ge {|A|}/{d^2}$.
\end{enumerate}

Let $\beta:={100d^4 \log  n }/{s}$. If $\beta^k<1/r$, then $G$ contains a cycle of length at least $\sum_{i\in [r]}d(B_i)$.
\end{proposition}
\begin{remark}\label{remark:reduction}  
In Proposition~\ref{lem:connecting}, if $\beta=1/s^{1-o(1)}$ and $k$ is a constant, then the condition $\beta^k<1/r$ can be replaced by the condition $s^k>n$. This follows from the observation that the disjointness of the sets $S^*_{i,j}$ implies that $r<n/s$.
\end{remark}
In Proposition~\ref{lem:connecting}, $\beta$ bounds the proportion of marked trees that remain after one round. After $k$ rounds, either a cycle has been formed or at most $\beta^k r<1$ marked trees remain. The latter is impossible, giving the required cycle.

The next theorem states the expansion properties of the $2$-core of the giant component of $Q_p^d$ that we will use in the application of Theorem~\ref{thm: cycle packing}.
\begin{theorem}\label{thm: expansion of 2-core}
Let $\eps>0$ be a sufficiently small constant and let $p=p(d) = (1+\eps)/d$.  Let \[\alpha=\alpha(\eps)\coloneqq \eps^{1/\eps^{20}}, \quad m=m(\eps,d) \coloneqq \eps^2\cdot  2^{\alpha d},  \quad \text{and} \quad M=M(\eps,d)\coloneqq   \eps^{-2}\cdot  2^{\alpha d}.\]
Let $C_2$ be the $2$-core of the giant component of $Q_p^d$. Whp every connected set $S\subseteq V(C_2)$ in $C_2$ with size 
$m\leq |S| \leq M+1$ satisfies the following:

\vbox{
\begin{enumerate}[leftmargin= \widthof{R100000}]
    \item $\eps|S|/4000 \leq |N_{C_2}(S)|\leq 2000 \eps |S|$, and 
    \item for every $Z\subseteq N_{C_2}(S)$ we have that $e_{C_2}(S,Z)\leq |Z|+ \eps^{8}\cdot  2^{\alpha d}$.
\end{enumerate}
}
\end{theorem}

\begin{proof}[Proof of Theorem~\ref{thm:main}]
Let $p_1=(1+2\eps/3)/d$. We couple $Q_{p_1}^d$ and $Q_p^d$ such that $Q_{p_1}^d \subseteq Q_p^d$. Let $\alpha=\alpha(\eps)\coloneqq \eps^{1/\eps^{20}}$ be as in the statement of Theorem~\ref{thm: expansion of 2-core}, and let $k\coloneqq\lceil 2/\alpha \rceil$, $h\coloneqq \lceil \log_2(k+1)\rceil$ and $D\coloneqq d-h-2$. Let $H_0\subseteq Q^d$ be the $(d-2)$-dimensional hypercube obtained by fixing the first $2$ coordinates to $1,1$. Let $H_1,\dots,H_k\subseteq Q^d$ be vertex-disjoint hypercubes of dimension $D$, obtained by fixing the first $h+2$ coordinates, with the first $2$ fixed to $0,0$. Order them so that, for $1\leq j <k$, $H_{j}$ is connected to $H_{j+1}$ by a perfect matching $M_j$ in $Q^d$. In addition, the fixed coordinates in $H_1$ are all $0$. Let $G_0\coloneqq Q_p^d[V(H_0)]$, and for $j\in [k]$ let $G_j^-\coloneqq Q_{p_1}^d[V(H_j)]$ and $G_j\coloneqq Q_{p}^d[V(H_j)]$, thus $G_j^-\subseteq G_j$.

\paragraph{Finding the paths in the $2$-core of $G_0$.}
Let $\cE$ be the event that the following are satisfied.
\begin{itemize}
\item[(P1)] The $2$-core of $G_0$  has a component  of order $(2\eps^2+O(\eps^3))2^{d-2}$, denoted by $C_2$, that satisfies the expansion properties described in Theorem~\ref{thm: expansion of 2-core}.
\item[(P2)] For $j\in [k]$, $G_j^-$ has a component $L_j^-$ such that every vertex in $V(H_j)$ is within distance at most $2$ from $V(L_j^-)$ in $H_j$.
\end{itemize}
Theorems~\ref{thm: size of 2-core},~\ref{thm: expansion of 2-core} and Lemma~\ref{lem: giant is well spread} imply that $\cE$ holds whp. We reveal $G_0$ and the graphs $G_j^-$, for $j\in [k]$, and assume that $\cE$ occurs. Observe that adding a vertex to an almost connected subset of $Q^d$ so that it becomes connected changes its vertex and edge boundaries by at most $d$. Thus (P1) gives that every almost connected set $S\subseteq V(C_2)$ in $C_2$ with size $m\leq |S| \leq M$ satisfies the following:

\vbox{
\begin{itemize}
    \item $|N_{C_2}(S)|\geq \eps|S|/4000-d\geq \eps |S|/5000$,
    \item $|N_{C_2}(S)| \leq 2000 \eps (|S|+1)+1\leq 2001\eps |S|$, and   
    \item for every $Z\subseteq N_{C_2}(S)$  we have $e_{C_2}(S,Z)\leq |Z|+\epsilon^{4}2^{\alpha d}+d$.
\end{itemize}
}
Applying \Cref{thm: cycle packing} with $k=2^{\alpha d}$, $c=\eps/5000$ and $C=2001\eps$ and then removing an edge from each cycle gives that, for some constant $A>0$, there is a family of vertex-disjoint paths, each of length at least $\eps^2 2^{\alpha d}$, spanning in total at least $A\eps^2 2^d$ vertices. Choose such a family $P_1,\dots P_r$ in $G_0$, and orient each path.

\paragraph{Choosing terminals and joining the paths.}  Let $q=1/(d\log d)^2$.
For every $j\in [k]$, every edge $e\in E(H_j)$ and every vertex $v\in V(H_j)$ within distance at most one from an endpoint of $e$, generate a variable $Y_{e,v}\sim\text{Bernoulli}(q)$. All these variables are mutually independent and independent of $Q_{p_1}^d$ and $G_0,G_1^-,\dots G_k^-$. There are $2D\leq 2d$ eligible vertices for each edge. For $e\in  \bigcup_{j=1}^kE(H_j)$, add $e$ to $\bigcup_{i\in [k]} G_i^-$ if at least one of its $2D$ associated variables is $1$, and let $G$ be the resulting graph. The edges of $G$ are independent, each present with probability $p_G\leq p_1+2Dq<p$. Thus we can couple these graphs so that $G_j^-\subseteq G[V (H_j )] \subseteq G_j \subseteq Q_{p}^d[V (H_j )]$ for every $j\in [k]$.

For $j\in [k]$ and $v\in V(H_j)$ let $P(v)$ be a path of length at most $2$ from $v$ to $L_j^-$ in $H_j$. Such a path exists in the event $\cE$ by (P2). For $v\in V(H_0)$, let $v_j$ be its projection to $H_j$, and let $Y_v\sim \operatorname{Bernoulli}(1/d)$, independently for all such $v$ and of all other variables. Let $P(v)$ be the path from $v$ to $v_1$ obtained by changing its non-zero coordinates among the first $h+2$ coordinates from $1$ to $0$, changing its second coordinate first and its first coordinate last. Thus all intermediate vertices start with $(1,0)$. Call $v$ successful if the following are satisfied.
\begin{enumerate}
    \item[(P1)] $E(P(v))\subseteq E(Q_p^d)$, and all of the edges of the form $v_iv_{i+1}$ are retained in $Q_p^d$ for $1\leq i< k$,
    \item[(P2)] $Y_{e,v_j}=1$ for $j\in[k]$ and $e\in E(P(v_j))$, and
    \item[(P3)] $Y_v=1$ and there does not exist a vertex $u\neq v$ such that $u$ and $v$ differ only on their first $h+2$ coordinates with $Y_u=1$.  
\end{enumerate} 
Then, each vertex $v\in V(H_0)$ is successful with probability at least $p^{h+2+(k-1)}\cdot q^{2k} \cdot d^{-1} (1-d^{-1})^{2^{h+2}}\geq d^{-8k}$. In addition, each edge and each Bernoulli random variable associated with an edge or vertex can affect the success of at most $2^{h+2}$ such vertices. McDiarmid's inequality (see Theorem~\ref{thm:mcdiarmid}) gives that whp for $1\leq i \leq r$ there exist sets $S_{i,0}^+$ and $S_{i,0}^-$ of successful vertices, each of size $ A \eps^2 2^{\alpha d} \cdot d^{-8k}/10 $ chosen from the first and last $|V(P_i)|/3$ vertices of $P_i$, respectively.

Let $\cE^+$ be the event that the following are satisfied.
\begin{enumerate}
    \item[(Q1)] The sets $S_{i,0}^+$, $S_{i,0}^-$, described above, exist for $i\in [r]$.
    \item[(Q2)] For $j\in [k]$, the second largest component of  $G_j$ has at most $d^2$ vertices and its largest component, denoted by $L_j$, satisfies the following property. For every set $S\subseteq V(L_j)$ of size $|S|\leq 2|V(L_j)|/3$ we have that $|N_{L_j}(S)|\geq |S|/d^2$.   
\end{enumerate} 
The preceding argument and Lemma~\ref{lem: giant weak expansion} imply that $\cE^+$ occurs whp. We reveal the remaining randomness and assume that $\cE^+$ occurs. For $i\in [r]$, $j\in [k]$,  and $*\in \{+,-\}$ let $S_{i,j}^*$ be the projection of  $S_{i,0}^*$ in $G_j$. In the event $\cE^+$, (P2) implies that the set $S_{i,j}^*$ lies in the component of $G_j$ containing $L_j^-$. Since $|S_{i,j}^*|> d^2$ and, by (Q2), the second largest component of $G_j$ has order at most $d^2$ it follows that $S_{i,j}^*$ lies in $L_j$.
In addition, by (P1), for $2\leq j\leq k$, $S_{i,j}^*$ is connected to $S_{i,j-1}^*$ via a matching in $Q_p^d$. For $1\leq i\leq r$ take $T_i$ to be the tree obtained by taking the union of the path $P_i$ and the paths $P(v)$ for $v\in S_{i,0}^+\cup S_{i,0}^-$. (P1) implies that $T_i\subseteq Q_p^d$. Then, (P3) implies that the trees $T_i$, $i\in [r]$ are vertex disjoint, as are the projected terminal sets in each
subcube. This gives a set of vertex disjoint marked trees $B_i^0=(T_i,S_{i,1}^+,S_{i,1}^-)$ each of length at least $A\varepsilon^2 2^{\alpha d}/3$. Let $L_1^+$ be the union of $L_1$ and these marked trees.

Applying Proposition~\ref{lem:connecting} (see also Remark~\ref{remark:reduction}) with $k-1$ rounds to $L_1^+$ and the giants $L_2,\dots,L_k$ with $s=  2^{\alpha d-O_\varepsilon(\log d)}$ and $k=\lceil 2/\alpha \rceil$  implies that $Q_p^d$ spans a cycle of length at least $A\varepsilon^2 2^d/3$, completing the proof. 
\end{proof}

\section{Proof of Theorem~\ref{thm: cycle packing}}\label{sec: path packing}
We use the following variant of Krivelevich's result on long cycles in locally expanding graphs (Theorem 1 of \cite{krivelevich2019long}).

\begin{proposition}\label{thm: krivelevich}
Let $h\geq 1,$ and $r\geq 2$. Let $G$ be a connected graph on more than $h$
vertices satisfying:
\[\text{$|N(W )| \geq r$, for every  almost connected $W\subseteq V(G)$ with $h/2\leq |W| \leq h$}.\]
Then $G$ contains a cycle of length at least $r+1$.    
\end{proposition}

\begin{proof}
Fix $v\in V(G)$, and let $T_v$ be a depth-first search (DFS) tree of $G$ rooted at $v$. By Proposition~2.1 of~\cite{krivelevich2019long}, there exist a vertex $u$ and a set $X$ of children of $u$ in $T_v$ such that, if $W$ denotes the union of the subtrees rooted at the vertices of $X$, then $h/2\le |W|\le h$. The set $W$ is almost connected, with witness $u$. Let $P_{v,u}$ be the $v$--$u$ path in $T_v$. By the ancestor--descendant property of DFS (Proposition~2.2 of~\cite{krivelevich2019long}), every external neighbour of $W$ lies on $P_{v,u}$; hence $N_G(W)\subseteq V(P_{v,u})$. By the hypothesis, $|N_G(W)|\ge r$. Let $a$ be the vertex of $N_G(W)$ closest to $v$ in $T_v$ and let $b\in W$ be adjacent to $a$. The $a$--$b$ path in $T_v$ contains all the vertices of $N_G(W)$ and at least one vertex of $W$. Adding the edge $ab$ therefore gives a cycle with at least $|N_G(W)|+1\geq r+1$ vertices.
\end{proof}

We now pass from a single cycle to a packing that covers a positive proportion of the graph, proving Theorem~\ref{thm: cycle packing}. For the proof outline, see Section~\ref{sec:outline}.

We now prove Theorem~\ref{thm: cycle packing}. For a sketch of its proof see \cref{sec:outline}.

\begin{proof}[Proof of  Theorem~\ref{thm: cycle packing}] 
Let $h:=\max\{k,ck\}$. Let $C_1, \dots ,C_t$ be a maximal family of vertex-disjoint cycles, each of length at least $\lceil \ell \rceil$ in $G$. Note that \ref{item: vertex expansion} implies that $|N(W)|\geq \ell$ for every  almost connected $W\subseteq V(G)$ with $k/2\leq |W| \leq k$. Thus, Proposition~\ref{thm: krivelevich} implies that $G$ has at least one cycle of length at least  $\lceil \ell \rceil$, so $t\geq 1$. For $i\in [t]$ let \[S_i\coloneqq V(C_i).\]
We now define pairwise disjoint subsets $S_{t+1},S_{t+2},...$ of $V(G)$ recursively. For $i\geq t+1$, given the sets $S_{t+1},\dots, S_{i-1}$, let $S_{\leq i-1}\coloneqq \cup_{j=1}^{i-1}S_j$ and $S_{i}$ be an almost connected subset  of $V(G)\setminus S_{\leq i-1}$ satisfying  
\begin{equation}\label{property S}
|S_{i}|\in [h/2,h],   \text{  and }\left|N(S_i)\setminus S_{\leq i-1} \right| < \ell.  \end{equation}
If no such set exists, then set $S_i=\emptyset$ and stop. We let $\tau$ be the largest index $i$ such that $S_{i}\neq \emptyset$.  Let
 \[
S\coloneqq S_{\leq \tau}, \hspace{5mm} S^-\coloneqq S_{\leq t}, \hspace{5mm} \text{ and }\hspace{5mm} U\coloneqq V(G)\setminus S.
 \]
By maximality of the cycle family $C_1,\dots C_t$,  $G[U]$ contains no cycle of length at least $\lceil \ell \rceil$. Since $S_{\tau+1}=\emptyset$, every almost connected set $X\subseteq U$ with $h/2\le |X|\le h$ satisfies  $|N_{G[U]}(X)|=|N_G(X)\setminus S|\ge \ell$. If a component $Q$ of $G[U]$ had more than $h$ vertices, then every set that is almost connected in $Q$ is also almost connected in $G$, and its  neighbourhood in $Q$ equals its  neighbourhood in $G[U]$. Hence,  Proposition~\ref{thm: krivelevich} yields a cycle of length at least $\lceil \ell \rceil$ in $Q$, contradicting the maximality of $C_1,\dots C_t$. Thus every component of $G[U]$ has at most $h$ vertices.

Partition the components of $G[U]$ according to their neighbourhoods in $G$, into $A_1,A_2,...,A_x$ and $ B_1,B_2,...,B_y$, where \[
|N_G(A_i)|\geq c|A_i|/2 \text{ for }i\in [x]
\hspace{5mm}\text{ and } \hspace{5mm}
|N_G(B_i)|< c|B_i|/2 \text{ for }i\in [y].
\]
Assign each poorly expanding component to the last removed set it meets: for $j\in[y]$ let $\sigma(j)$ be the maximum index $i\in[\tau]$ such that $B_j$ is incident to $S_i$. Such an index exists since $G$ is connected and $S\neq\varnothing$. For $i\in[\tau]$ let $T_i:=S_i\cup\bigcup_{j:\,\sigma(j)=i}B_j$, and let
\[
A:=\bigcup_{i=1}^{x}A_i,\qquad
T:=\bigcup_{i=1}^{\tau}T_i,\qquad
T^-:=\bigcup_{i=1}^{t}T_i.
\]
Then $V(G)=T^-\uplus(T\setminus T^-)\uplus A$.
\begin{claim}
The following inequalities hold.
\begin{equation}\label{eq: bound T-}
    |T^-| \leq \left(\frac{2C}{c}\right)|S^-|.
\end{equation}

\begin{equation}\label{eq: bound S minus S- and A}
    |T\setminus T^-|+|A| \leq \left(\frac{5C}{c}\right)|S^-|.
\end{equation}
\end{claim}

\begin{proof}
First we prove that $|T_i|\leq 2(C/c)|S_i|$ for $i\in[\tau]\setminus[t]$. Indeed, assume otherwise. Since every $B_j$ has at most $h$ vertices, there is a subset $\phi(i)\subseteq\sigma^{-1}(i)$ such that $T_i':=S_i\cup\bigcup_{j\in\phi(i)}B_j$ has size in  $[2(C/c)|S_i|,2(C/c)|S_i|+h]\subseteq [(C/c)h,2(C/c)h+h]\subseteq [(C/c)h,3(C/c)h]\subseteq I$. The set $T_i'$ is almost connected, since each selected component is incident to $S_i$. Thus \ref{item: vertex expansion} gives $|N_G(T_i')|\geq c|T_i'|$. But
\begin{align*}
    |N_G(T_i')|& \leq |N_G(S_i)|+ \sum_{j\in \phi(i)}|N_G(B_j)| < C|S_i|+ \frac{c|T_i'|}2  \leq c|T_i'|,
\end{align*}
a contradiction.

The same argument applies to any connected segment $R$ of a cycle with $ck/10\leq |R|\leq h$, together with any of the $B_j$ assigned to that cycle and incident to $R$. Partition each $C_i$, $i\in[t]$, into consecutive segments with between $ck/10$ and $ck$ vertices. This is possible since $|C_i|\geq ck/10$. Assign each $B_j$ with $\sigma(j)=i$ to one segment incident to it. Applying the preceding argument to each segment and its assigned components, and summing, gives $|T_i|\leq 2(C/c)|S_i|$ also for $i\in[t]$. Summing over $i\in[t]$ gives~\eqref{eq: bound T-}.

To prove~\eqref{eq: bound S minus S- and A}, let $m$ be the number of edges from $S^-$ to $V(G)\setminus S^-$. First, for each of the cycle segments $R$ above, \ref{item: vertex expansion} and  \ref{item: vertex implies edge} give that $e_G(R,V(G)\setminus R)\leq C|R|+c^2k/10\leq C|R|+c|R|\leq 2C|R| $. Consequently, summing over the segments gives
\begin{equation}\label{eq:initial-edges}
m\leq 2C|S^-|.
\end{equation}

For $i>t$, let $T_{<i}\coloneqq \bigcup_{j<i}T_j$, $T_{>i} \coloneqq \bigcup_{j>i}T_j$ and $S_{<i}\coloneqq S_{\leq i-1}$. Counting edges between earlier and later parts gives 
\[m\geq e_G(S^-,(T \setminus T^
-) \cup A) \geq  \sum_{i=t+1}^\tau e_G(T_i,T_{<i})+e_G(A,S)- \sum_{i=t+1}^\tau e_G(T_i,T_{>i}\cup A).\]
Now, since $T_i$ is almost connected and $ck/10\leq |T_i|\leq 2C|S_i|/c\leq 2Ch/c$, \ref{item: vertex expansion} applies to $T_i$. By the definition of $\sigma$, every neighbour of a component assigned to $T_i$ lies in $S_{\leq i}$. Thus~\eqref{property S} gives
\[ e_G(T_i,T_{<i})\geq |N_G(T_i)\cap T_{<i}|\geq c|T_i|-|N(S_i)\setminus S_{<i}|\geq  c|T_i|-\frac{ck}{10} \geq \frac{4c|T_i|}{5}. \]
Also, 
\ref{item: vertex implies edge} and \eqref{property S} give
\[
e_G(T_i,T_{>i}\cup A)
=e_G(S_i,T_{>i}\cup A)
\leq e_G(S_i,N(S_i)\setminus S_{<i})\leq \frac{ck}{10}+\frac{c^2k}{10}\leq \frac{2c|S_i|}5\leq \frac{ 2c|T_i|}5.
\]
Finally, the $A_i$ are components of $G[U]$, so
\[
e_G(A,S)=\sum_{i=1}^{x}e_G(A_i,S)
\geq\sum_{i=1}^{x}|N_G(A_i)|\geq\frac{c}{2}|A|.
\]
Combining the above inequalities gives
\begin{align*}
2C|S^-|&\geq \sum_{i=t+1}^\tau e_G(T_i,T_{<i})+e_G(A,S)- \sum_{i=t+1}^\tau e_G(T_i,T_{>i}\cup A)
\\&\geq \sum_{i=t+1}^\tau \frac{4c|T_i|}5+\frac{c|A|}2- \sum_{i=t+1}^\tau \frac{2c|T_i|}5 \geq \frac{2c(|A|+|T\setminus T^-|)}5. 
\end{align*}
Rearranging the above inequality gives \eqref{eq: bound S minus S- and A}.
\end{proof}
\eqref{eq: bound T-} and \eqref{eq: bound S minus S- and A} give that $|V(G)|=|T^-|+|T\setminus T^-|+|A|\leq 7C|S^-|/c$. Thus $C_1,...,C_t$ cover at least $(1/7)(c/C)|V(G)|$ vertices, completing the proof. 
\end{proof}

\section{Proof of Proposition~\ref{lem:connecting}}
For a sketch of the proof of Proposition~\ref{lem:connecting}, see Section~\ref{sec:outline}. We first analyze one round, comparing the number of vertices used by the connecting paths with the expansion available to the unused terminal sets.

\begin{lemma}\label{lem:one-round} Let $n,d,s \geq 1$, and $r\geq 2$ be integers. Let $G$ be an $n$-vertex graph, and let $S_1,\dots,S_{2r}\subseteq V(G)$ be pairwise disjoint sets of size $s$.  Assume that
\begin{equation}
    |N_G(A)|\ge \frac{|A|}{d^2}
  \qquad\text{for every }A\subseteq V(G)\text{ with }|A|\le 2n/3.
  \label{eq:expansion}
\end{equation}
Let $M$ be a perfect matching on the set $[2r]$. Let $\cM^+$ be the set of matchings $M^+$ on $[2r]$ satisfying the following conditions.

\begin{enumerate}[label = \textnormal{(\Roman*)}, leftmargin= \widthof{R100000}]
\item $M^+$ is a matching on $[2r]$ that is edge disjoint from $M$, and the graph $([2r],M\cup M^+)$ is either a Hamilton cycle or acyclic.\label{item: acyclic}
\item There is a family of vertex-disjoint paths $\cP=\{P_{xy}:xy\in M^+\}$ such that,
for every $xy\in M^+$, the path $P_{xy}$
is an $S_x$-$S_y$ path in $G$ of length at most
      \[
        L:=20d^2\log n-1.
      \]
In addition, with $U \coloneqq  \bigcup_{P\in \cP} V(P)$, we have that $U \cap S_z = \emptyset$ for every $z  \in [2r] \setminus  V(M^+)$, i.e. the paths in $\cP$ avoid every set $S_z$ such that $z$ is not matched by $M^+$. 
\label{item: path family}
\end{enumerate}
Let $M^+\in \cM^+$ be of maximum cardinality. If $([2r],M\cup M^+)$ is not a Hamilton cycle, then
\begin{equation}\label{eq:one-round-bound}
   r-|M^+|\le \frac{100d^4 \log n}{s}\,r.   
\end{equation}
 \end{lemma}

\begin{proof}[Proof of Lemma~\ref{lem:one-round}] Let $M^+\in\cM^+$ have maximum cardinality, and let
\[
    \Gamma:=([2r],M\cup M^+),
    \qquad m:=|M^+|,
    \qquad q:=r-m.
\]
If $\Gamma$ is a Hamilton cycle, then $|M^+|=r$ and \eqref{eq:one-round-bound} holds trivially. Otherwise, condition~\ref{item: acyclic} implies that $\Gamma$ is acyclic. Since $M$ is a perfect matching and $M^+$ is a matching, every vertex in $\Gamma$ has degree~$1$  or $2$. Hence every component of $\Gamma$ is a path of length at least~$1$. Consequently, the number of components of $\Gamma$ is
$|V(\Gamma)|-|E(\Gamma)|=2r-(r+m)=q.$

Let $X$ be the set of endpoints of the components of $\Gamma$. Since $M$ is a perfect matching, these endpoints are precisely the vertices not covered by $M^+$; that is, $X=[2r]\setminus V(M^+).$ In particular, $|X|=2q\ge 2$. We partition these endpoints so that an edge across the two parts either joins distinct path components or closes the final cycle. If $q=1$, then $\Gamma$ is a Hamilton path, and we let $X^+$ and $X^-$ be singleton sets containing its two distinct endpoints. If $q\ge 2$, let $X^+$ be the set of endpoints of $\lfloor q/2\rfloor$ components of $\Gamma$, and let $X^-:=X\setminus X^+$. In both cases, $X=X^+\uplus X^-$, with $|X^+|,|X^-|\ge \frac{2q}{3}$, and, for every $x^+\in X^+$ and $x^-\in X^-$, the graph $\Gamma+\{x^+x^-\}$ is either a Hamilton cycle or acyclic.

Since $M^+\in\cM^+$, there exists a family $\cP=\{P_{xy}:xy\in M^+\}$ of vertex-disjoint paths of length at most $L$ such that, for every $xy\in M^+$, the path $P_{xy}\subseteq G$ has one endpoint in $S_x$ and the other in $S_y$, and its vertex set avoids $S_z$ for every $z\in[2r]\setminus V(M^+)$. Let $    U:=\bigcup_{xy\in M^+}V(P_{xy})$. Then
\begin{equation}\label{eq: bound U}
    |U|\le m(L+1)=20d^2m\log n\le 20d^2r\log n, 
\end{equation}
where the last inequality follows from $m=|M^+|\le r$.

For $\ast\in\{+,-\}$, let $S(X^\ast):=\bigcup_{i\in X^\ast}S_i$. Also let $S(X)=S(X^+)\cup S(X^-)$. Since $X=[2r]\setminus V(M^+)$, the choice of $\cP$ gives
\begin{equation}\label{eq: properties of S(X) }
    |S(X^\ast)|\ge \frac{2q}{3}s
    \qquad\text{and}\qquad
    S(X^\ast)\cap U=\varnothing.  
\end{equation}
Suppose, for a contradiction, that \eqref{eq:one-round-bound} fails; that is,
$ q=r-m>r(100d^4\log n)/s.$ We show that the remaining terminal sets can still be joined by a short path avoiding $U$. First, by the expansion hypothesis and \eqref{eq: bound U}, for every $W\subseteq V(G)\setminus U$ with $2qs/3\le |W|\le 2n/3$, we have
\begin{equation}\label{eq: G-U expansion} 
    |N_{G-U}(W)|
    =|N_G(W)\setminus U| \ge \frac{|W|}{d^2}-20d^2r\log n \ge \frac{|W|}{d^2}-\frac{qs}{5d^2}
     \ge \frac{|W|}{2d^2}.
\end{equation}

We now grow neighbourhoods in $G-U$ from the two groups of terminals. Let $S_0^-:=S(X^-)$ and, for $i\ge 1$, define $S_i^-:=N_{G-U}(S_{i-1}^-)\cup S_{i-1}^-$.
Since $2qs/3\le |S_0^-|\le |S_{i-1}^-|$, \eqref{eq: G-U expansion} implies that, if $|S_{i-1}^-|<2n/3$, then
\[
    |S_i^-|
    =|N_{G-U}(S_{i-1}^-)|+|S_{i-1}^-|
    \ge \left(1+\frac{1}{2d^2}\right)|S_{i-1}^-|.
\]
It follows by induction that, with $t:=\lfloor L/2\rfloor$,
\begin{equation*}
    |S_t^-|\ge
    \min\left\{
        \left(1+\frac{1}{2d^2}\right)^t|S_0^-|,
        \frac{2n}{3}
    \right\} \ge
    \min\left\{
        \left(1+\frac{1}{2d^2}\right)^t,
        \frac{2n}{3}
    \right\}
    =\frac{2n}{3}. 
\end{equation*}
Indeed, since $0\le 1/(2d^2)\le 1$, the inequality
$\log(1+x)\ge x-x^2/2\ge x/2$ gives
\[
    \left(1+\frac{1}{2d^2}\right)^t
    =\exp\left(t\log\left(1+\frac{1}{2d^2}\right)\right)
    \ge \exp\left(\frac{t}{4d^2}\right)
    \ge \exp(\log n)=n.
\]
For the penultimate inequality, note that $n\geq 2rs\geq 4$ and hence $t=\lfloor L/2\rfloor\geq 10d^2\log n- 3/2 \ge 4d^2\log n$.

Defining $S_i^+$ analogously gives $|S_t^+|\ge 2n/3$. Hence
$S_t^+\cap S_t^-\ne\varnothing$. Let $v\in S_t^+\cap S_t^-$. Then there is an $S(X^+)$--$v$ path of length at most $t$ and a $v$--$S(X^-)$ path of length at most $t$, both in $G-U$. Their concatenation contains an $S(X^+)$--$S(X^-)$ path $P$ of length at most $2t\le L$.

Let $P'$ be a shortest subpath of $P$ with one endpoint in $S(X^+)$ and the other in $S(X^-)$. Denote its endpoints by $u^+$ and $u^-$, and let $x^+\in X^+$ and $x^-\in X^-$ be the unique indices such that $u^+\in S_{x^+}$ and $u^-\in S_{x^-}$. By the minimality of $P'$, its vertex set is disjoint from $S(X)\setminus\{u^+,u^-\}$. Indeed, if an internal vertex of $P'$ belonged to $S(X^+)$ or $S(X^-)$, then an appropriate proper subpath of $P'$ would be a shorter $S(X^+)$--$S(X^-)$ path.

By construction, $\Gamma+\{x^+x^-\}$ is either a Hamilton cycle or acyclic. The edge $x^+x^-$ is not in $M$: if $q\ge 2$, then $x^+$ and $x^-$ lie in distinct components of $\Gamma$, while if $q=1$, they are the endpoints of a Hamilton path on $2r\ge 4$ vertices. Moreover, $x^+$ and $x^-$ are not covered by $M^+$, so $M^+\cup\{x^+x^-\}$ is a matching edge-disjoint from $M$. Thus condition~\ref{item: acyclic} holds. Since $P'\subseteq G-U$, the set $\cP\cup\{P'\}$ is a family of vertex-disjoint paths. Furthermore, both $U$ and $V(P')$ are disjoint from  $S_z$ for every
$z\in X\setminus\{x^+,x^-\}$, so this enlarged family also satisfies condition~\ref{item: path family}. Therefore, $M^+\cup\{x^+x^-\}\in\cM^+$,
contradicting the maximality of $M^+$. Hence  \eqref{eq:one-round-bound} holds.
\end{proof}

\begin{proof}[Proof of Proposition~\ref{lem:connecting}] Let $\ell=\sum_{i\in [r]}d(B_i)$. 
Re-index the sets so that $S_{i,j}=S_{i,j}^+$ and  $S_{i+r,j}=S_{i,j}^-$ for $i\in[r]$ and $j\in\{0,1,\ldots,k\}$, and re-index the matchings in the same way. We prove by induction on $j\in\{0,1,\ldots,k\}$ the following statement,
denoted by $P(j)$. Either $G$ contains a cycle of length at least $\ell$, or
there is a family of at most $\beta^j r$ pairwise disjoint marked trees in $G$, with all vertices in $\cup_{i=0}^jV(G_i)$, and with
terminal sets in $\{S_{i,j}:i\in[2r]\}$ whose total length is at least
$\ell$.

The base case $P(0)$ follows from condition~\textup{(a)} by taking $ \cB_0=\{(T_i,S_{i,0},S_{i+r,0}):i\in[r]\}$.

Let $0\le j<k$ and assume that $P(j)$ holds. If its first alternative holds,
there is nothing to prove. Otherwise, let $\cB_j=\{B_i=(T_i,S_{a_i,j},S_{b_i,j}):i\in[t]\}$.

We first handle the case $t=1$. Condition~\textup{(c)} implies that
$G_{j+1}$ is connected. Let $Q$ be a shortest
$S_{a_1,j+1}$--$S_{b_1,j+1}$ path in $G_{j+1}$, with endpoints $u\in S_{a_1,j+1}$ and $v\in S_{b_1,j+1}$.
Let $u'\in S_{a_1,j}$ and $v'\in S_{b_1,j}$ be matched to $u$ and $v$ by
$M_{a_1,j}$ and $M_{b_1,j}$, respectively. Then $Q$, the edges $uu'$ and
$vv'$, and the unique $u'$--$v'$ path in $T_1$ form a cycle of length at
least $d(B_1)\ge \ell$. Thus $P(j+1)$ holds.

We may therefore assume that $t\ge2$. Relabel the $2t$ terminal sets and the associated matchings by
\[
  A_i:=S_{a_i,j+1},\qquad A_{i+t}:=S_{b_i,j+1},\qquad
  F_i:=M_{a_i,j},\qquad F_{i+t}:=M_{b_i,j}
  \qquad(i\in[t]),
\]
and let $M:=\{\{i,i+t\}:i\in[t]\}$. Apply Lemma~\ref{lem:one-round} in $G_{j+1}$ to
$A_1,\ldots,A_{2t}$ and $M$. Let $M^+$ have maximum cardinality among the matchings allowed by that lemma, let $\cP$ be a witnessing family of paths, and put $\Gamma:=([2t],M\cup M^+)$. If $\Gamma$ is acyclic, then, writing $t':=t-|M^+|$, Lemma~\ref{lem:one-round} and $|V(G_{j+1})|\leq n$ give $t'\le \beta t.$

The graph $\Gamma$ records how the current marked trees are joined. For every $x\in V(M^+)$, let $u_x\in A_x$ be the endpoint in $A_x$ of the unique path in $\cP$ incident with $x$, and let $v_x$ be the vertex matched to $u_x$ by $F_x$.

If $\Gamma$ is a Hamilton cycle, then the paths in $\cP$, the matching edges $u_xv_x$, and, for each $i\in[t]$, the unique $v_i$--$v_{i+t}$ path in $T_i$ form a cycle. Its length is at least $\sum_{i=1}^t d(B_i)\ge \ell,$ so $P(j+1)$ holds.

Suppose instead that $\Gamma$ is acyclic. Its $t'$ components are paths,
each with two endpoints not covered by $M^+$. For each component $C$, form
a graph $T_C$ by taking the trees $T_i$ corresponding to the edges of $M$
in $C$, the paths in $\cP$ corresponding to the edges of $M^+$ in $C$, and
the matching edges $u_xv_x$ for $x\in V(C)\cap V(M^+)$. At each of the two
endpoints $x$ of $C$, also add the entire perfect matching $F_x$. By condition~\textup{(II)} of
Lemma~\ref{lem:one-round}, the corresponding sets $A_x$ are disjoint from all paths in $\cP$. Since $C$ is a
path and all constituent trees and paths are otherwise vertex-disjoint, $T_C$ is a tree; moreover, the trees obtained from
distinct components are vertex-disjoint.

Taking the two sets $A_x$ indexed by the endpoints $x$ of $C$ as the terminal sets of $T_C$ gives a marked tree. Every path between these terminal sets traverses each marked tree of $\cB_j$ corresponding to an edge of $M$ in $C$, from one of its terminal sets to the other. Consequently, the total length of the resulting $t'$
marked trees is at least the total length of $\cB_j$, and hence at least $\ell$.
Moreover, $t'\le \beta t\le \beta^{j+1}r.$
Thus $P(j+1)$ holds, completing the induction.

Finally, apply $P(k)$. If $G$ contained no cycle of length at least $\ell$,
there would be a family of at most $\beta^k r<1$ marked trees whose total length
is at least $\ell$. Such a family would be empty, a contradiction. Hence
$G$ contains a cycle of length at least $\ell$.
\end{proof}

\section{Properties of the 2-core}\label{sec: proof expanstion}
Throughout this section, let $\eps>0$ be a sufficiently small constant and let $p=(1+\eps)/d$. Let $L_1$ denote the largest component of $Q_p^d$ and let $C_2$ denote the largest component of its $2$-core.

This section has two goals. We first estimate $|V(C_2)|$, proving Theorem~\ref{thm: size of 2-core}. The proof uses local explorations to identify configurations that certify whether a vertex can belong to the core. We then prove Theorem~\ref{thm: expansion of 2-core}. For this, we consider spanning trees of connected sets in $C_2$ and attach their boundary vertices as leaves. Explorations from the leaves and from the interior of such a tree have different constraints: the leaves need further connections in the core, whereas edges from an interior vertex to vertices outside the tree can only lead outside the core.

\subsection{The order of the 2-core}
We begin with the heuristic behind Theorem~\ref{thm: size of 2-core}. A vertex $v$ should belong to $C_2$ when it has two neighbours in $Q_p^d$ whose explorations in $Q_p^d-\{v\}$ survive for long enough. Each exploration has survival probability approximately $2\eps$, and the explorations are approximately independent. This suggests that the probability that $v$ belongs to $C_2$ is $\Pr(\operatorname{Bin}(d,2\eps p)\geq2)\approx2\eps^2$, and hence that $|V(C_2)|$ is approximately $2\eps^2 2^d$.

The proof makes this heuristic precise in two steps. For the upper bound, a vertex in the core must either belong to a small connected subgraph containing a cycle or have two neighbours whose explorations reach a fixed large size. Estimating the probability of these local configurations gives the desired upper bound for a fixed vertex. For the lower bound, we construct two disjoint explorations from neighbours of $v$, each of which reaches many vertices lying on cycles. These configurations place $v$ in a large component of the $2$-core and therefore in $C_2$ whp. In both directions, we use McDiarmid's inequality to establish concentration of the number of vertices satisfying the corresponding properties.

\begin{proof}[Proof of Theorem~\ref{thm: size of 2-core}]
We first prove the upper bound. Let $r=\eps^{-6}$ and call a vertex $v$ good if it lies in a connected subgraph of $Q_p^d$ of order at most $2r+1$ containing a cycle, or if at least two of its neighbours in $Q_p^d$ have components of order at least $r$ in $Q_p^d-\{v\}$. Every vertex of $C_2$ is good.

For a fixed pair $u,u'$ of neighbours, explore their components in $Q_p^d-\{v\}$ successively, stopping each exploration when it reaches $r$ vertices and never querying an edge to a previously reached vertex. If $v$ is good and does not lie in a connected subgraph of order at most $2r+1$ containing a cycle, then both explorations reach $r$ vertices for some such pair. By \eqref{eq: exploring to medium}, the probability that both explorations reach this size is $4\eps^2+O(\eps^3)+o_d(1)$. Thus the probability that $v$ is good is at most
\[
\sum_{i=3}^{2\eps^{-6}+1} (ed)^{i-1}p^{i-1} \binom{i}{2}p+ \binom{d}{2}p^2(4\eps^2+O(\eps^3))=
2\eps^2+O(\eps^3).\]
An application of McDiarmid's inequality gives that whp there are at most
$(2\eps^2+O(\eps^3))2^{d}$ good vertices, hence  $|V(C_2)|\leq (2\eps^2+O(\eps^3))2^{d}$ whp.

We now prove the lower bound. Fix $v\in V(Q^d)$ and let $L=\log^2(d-2)$. Let $Q^d(v)$ be the $(d-2)$-dimensional subcube induced by the vertices that agree with $v$ in their first two coordinates, and let $Q_p^d(v)=Q_p^d[V(Q^d(v))]$. 

Expose the edges incident to $v$ in $Q_p^d(v)$. Form $N_v$ by choosing $L$ of its neighbours, or all of them if there are fewer than $L$. Let $R(v)=\{i\in [d]: \exists u\in N_v$ such that $u$ and $v$ differ in their $i$th coordinate$\}$ be the set of coordinates on which $v$ differs from its neighbours in  $N_v$. For $u\in N_v$, let $Q(u,v)$ be the percolated subcube obtained from $Q_p^d(v)$ by fixing the coordinates in $R(v)$ to their values at $u$. Its underlying cube has dimension $d-|R(v)|-2$. Let $N_{u,v}$ be the set of vertices within distance $L$ from $u$ in $Q(u,v)$. Let $\cE(u,v)$ be the event that $|N_{u,v}|\geq d^{20}/2$ and at least $d^2$ vertices $w\in N_{u,v}$ lie on the $4$-cycle $w_{0,0}w_{0,1}w_{1,1}w_{1,0}w_{0,0}$ in $Q_p^d$, where $w_{i_1,i_2}$ is obtained from $w$ by changing coordinate $j\in\{1,2\}$ whenever $i_j=1$. Let $\cE_v$ be the event that there exists a pair of vertices $u,u'\in N_v$ such that the events  $\cE(u,v)$ and $\cE(u',v)$ occur. Note that in the event $\cE_v$ the vertex $v$ lies in a component of the $2$-core of $Q_p^d$ with at least $2d^2$ vertices.   

By \eqref{eq: expansion}, for $u\in N_v$ the probability that at least $d^{20}$ vertices lie at distance at most $L$ from $u$ in $Q(u,v)$ is $x=2\eps+O(\eps^2)$. Since the corresponding $4$-cycles are vertex disjoint, each of them belongs to $Q_{p}^d$ independently with probability $p^4$. Thus, by the Chernoff bound at least $d^{2}$ of these cycles belong to $Q_{p}^d$ whp. It follows that 
\begin{align*}
\Pr(\cE_v)  &\geq \binom{d-2}{2}(xp)^2(1-xp)^{d-4} - \Pr(Bin(d-2,p)\geq \log^2 d)+o(1)
\\&= \frac{(1+o(1))(xpd)^2e^{-xpd}}{2}
= \frac{((2\eps+O(\eps^2))(1+\eps))^2e^{-(2\eps+O(\eps^2))(1+\eps)}}{2}
\\&= (2\eps^2+O(\eps^3))e^{-2\eps+O(\eps^2)}=
2\eps^2+O(\eps^3).
\end{align*}

The event $\mathcal E_v$ depends only on edge indicators within cube distance
$L+3$ of $v$. Changing one edge therefore affects  at most $2d^{L+3}\leq e^{2\log^3d}$ of those events. McDiarmid's inequality (see Theorem~\ref{thm:mcdiarmid}) gives that the event $\cE_v$ occurs for at least 
$(2\eps^2+O(\eps^3))2^d$ vertices $v\in V(Q_p^d)$ whp. Each of those vertices lies in a component of the $2$-core of $Q_p^d$ with at least $2d^2$ vertices, hence by Lemma~\ref{lem: giant weak expansion} in $C_2$ whp. This completes the proof of the lower bound on $|V(C_2)|$. 
\end{proof}

\subsection{Expansion properties of connected subsets of the 2-core}

We first discuss why connected sets in $C_2$ should expand by a factor of order $\eps$. The heuristic for Theorem~\ref{thm: size of 2-core} suggests that a vertex has degree at least $k$ in the core when at least $k$ of the explorations from its neighbours survive. Thus, for each fixed $k\geq3$, we expect only $O(\eps^k2^d)$ vertices of degree at least $k$, compared with $\Theta(\eps^2 2^d)$ vertices in the core altogether. If the degree-two vertices are spread roughly uniformly among the bare paths joining vertices of higher degree, then $C_2$ resembles an almost $3$-regular graph whose edges have been replaced by paths of typical length $\Theta(1/\eps)$. This suggests that expansion by a constant factor before these subdivisions becomes expansion by a factor of order $\eps$ in $C_2$.

This picture resembles the anatomy of the $2$-core of the binomial random graph, for which a much finer description was proved by Ding, Lubetzky and Peres~\cite{ding2014anatomy}.

A standard way to study the boundary of a connected set is to choose a spanning tree as a witness, count the possible trees, and bound the probability that a given witness has an unusually small or large boundary. Applied directly in $Q_p^d$, however, this argument does not capture the scale we need inside $C_2$: a tree may have a neighbourhood of order $|V(T)|$ in $Q_p^d$ even when its neighbourhood in the core has order $\eps|V(T)|$. We therefore need witnesses that also record the constraints imposed by membership in the core.

Start with a spanning tree of a connected set in $C_2$ and attach its neighbours in $C_2$ as leaves. The resulting tree $T'$ is \emph{interior-closed} in $C_2$: every neighbour in $C_2$ of a non-leaf vertex already belongs to $T'$. Explorations leaving the interior must therefore die out outside the core. By contrast, an exploration from a leaf must grow, find a cycle, or return to the tree through an additional edge. Heuristically, extinction at an interior vertex contributes a factor of about $1-2\eps$, while growth from a leaf contributes a factor of about $2\eps$. Thus, for a tree with $t'$ vertices and $s$ leaves, the expression $p^{t'-1}(1-2\eps)^{t'-s}(2\eps)^s$ suggests the probability scale for such a witness. The proof below handles the possible cycles, returns and dependencies through explicit explorations.

To obtain sufficiently strong bounds, we further restrict the form of the witness. For example, for the upper bound on expansion, we prune paths from an interior-closed tree while preserving its number of leaves. The resulting subtree has the following property: every leaf either has a neighbour of degree at least three in the tree or is incident to an additional edge with both endpoints in the tree.  We call such trees \emph{special}. Harper's inequality, together with a union bound, shows that connected sets in the relevant size range typically span few additional edges. Leaves of the second type are therefore rare. Leaves of the first type impose a useful restriction when counting trees: their incident edges are not subdivided when the tree is reconstructed from its contracted form. This gives the sharper count needed for the expansion bounds.

We first define special $(a,b)$-trees, where $a$ and $b$ count these two types of leaves. Lemma~\ref{lem: count prune trees} bounds their number, and Lemma~\ref{lem: special trees} shows how to prune or extend an interior-closed tree to obtain a special one. Then, Lemma~\ref{lem: existance of special trees} bounds the possible numbers of leaves of special trees in $C_2$. We finally combine these bounds to prove Theorem~\ref{thm: expansion of 2-core}.

\begin{definition} 
For a tree $T$ on at least three vertices, denote by $A(T)$ the set of leaves of $T$ whose unique neighbor in $T$ has degree at least three,  and by $B(T)$  the remaining leaves. 

For a graph $H\supseteq T$, call $T$ an
\emph{$H$-special $(a,b)$-tree} if $|A(T)|=a$, $|B(T)|=b$, and every leaf in $B(T)$
 is incident to an edge of $E(H[V(T)])\setminus E(T)$. The edges of $E(H[V(T)])\setminus E(T)$ are called \emph{$(H,T)$-surplus edges}.
 
We say that  $T$ is \emph{interior-closed in $H$} if $N_{H}(v)\subseteq V(T)$ for every non-leaf $v$ of $T$. It is \emph{almost interior-closed in $H$} if the same condition holds for every non-leaf that is not a neighbor of a leaf in $T$.
\end{definition}

For a tree $T$, write $\ell(T)$ for its number of leaves.
 
\begin{lemma}\label{lem: count prune trees}
Let $H$ be a graph of maximum degree at most $d$, let $v\in V(H)$, and let $k\ge3$ and $a,b\ge0$ be integers. Let $t(v,a,b,k)$ be the number of $H$-special $(a,b)$-trees rooted at $v$, with $k$ vertices, in which $v$  has degree at least $3$.
\begin{equation}\label{eq:count special trees}
 \text{If $k\ge 4(a+b)$, then } t(v,a,b,k) \leq (a+b) d^{k-1}\bfrac{2e^3k}{a+2b}^{a+2b}.  
\end{equation}
\end{lemma}

\begin{proof}
Let $\ell=a+b$ be the number of leaves of $T$. 
Fix an ordering of $V(H)$. For each vertex $u$, label its neighbours injectively by elements of $[d]$ according to this ordering. Direct the edges of a counted tree $T$ away from its root $v$, and label each directed edge $uw$ by the label of $w$ at $u$.

Contract every maximal bare subpath of $T$ into a single edge. Record on each resulting edge the sequence of labels along the corresponding path, and let $T'$ be the resulting rooted tree. It has $\ell$ leaves and at most $\ell-2$ vertices of degree at least three, including the root. Thus, with $h\coloneqq |V(T')|$, we have that
$\ell+1\le h \le 2\ell-1$.

The outgoing edges of $T'$ at each vertex have distinct first labels, inherited from the corresponding outgoing edges in $H$. Lemma~\ref{lem: tree count} gives that the number of such first-labeled rooted trees is at most $(ed)^{h-1}$.

Given $T'$, to determine $T$ first choose the $a$ leaves in $A(T')$ in at most $\binom{\ell}{a}\leq 2^\ell$ ways.  Their incident
edges have length one before and after contraction. The remaining length
assignments number
\[
 \binom{k-a-2}{h-a-2}
 \le\binom{k}{a+2b}
 \le\left(\frac{ek}{a+2b}\right)^{a+2b}.
\]
Here we used that $h-a-2\le a+2b$ and $a+2b\le k/2$. The remaining labels can be chosen in at most $d^{k-h}$ ways.
 Thus,
\[
t(v,a,b,k)
 \le \sum_{h=\ell}^{2\ell-1} (ed)^{h-1} 2^\ell \bfrac{ek}{a+2b}^{a+2b} d^{k-h}\leq 
  \ell d^{k-1}  \bfrac{2e^3k}{a+2b}^{a+2b}.
\]
\end{proof}

\begin{lemma}\label{lem: special trees}
Whp every tree $T\subseteq C_2$ with at least $d^3$ vertices that is interior-closed in $C_2$ admits trees
\[
T^-\subseteq T\subseteq T^+\subseteq C_2
\]
with the following properties.
\begin{enumerate}[label=(\roman*)]
\item $|V(T^-)|\ge |V(T)|/d^3$, $\ell(T^-)=\ell(T)$, and $T^-$ is $Q^d_p$-special.
\item $|V(T^+)|\le d^3|V(T)|$, $\ell(T)\le\ell(T^+)\le2\ell(T)$, and $T^+$ is $Q^d_p$-special and almost interior-closed in $C_2$.
\end{enumerate}
\end{lemma}

\begin{proof}
It suffices to construct $T^-$ and $T^+$ in the high-probability event given in Corollary~\ref{cor: sizes of paths bare paths}, that is, in the event that every component of $Q_p^d-V(C_2)$ has at most $d^3$ vertices and every bare path in $C_2$ has at most $d^3-3$ vertices.

To obtain $T^-$, let $F=E(Q_p^d[V(T)])\setminus E(T)$. Recursively delete a leaf whose parent has degree two in the current tree and which is not incident to $F$, until no such leaf exists. Let $T^-$ be the resulting tree. No endpoint of an edge in $F$ is deleted. Thus every remaining leaf with a degree-two parent is incident to a surplus edge whose endpoints both remain, so $T^-$ is special. At each step the parent becomes a leaf, so the number of leaves is unchanged.

This process cannot leave fewer than three vertices. Indeed, since $T$ has more than $d^3-1$ vertices, it is not a bare path in $C_2$. Thus it has a non-leaf vertex $v$ of degree at least three in $C_2$. Since $T$ is interior-closed, the neighbors of $v$ lie in $V(T)$. Either $v$ has degree at least three in $T$, in which case none of its neighbors is deleted, or $v$ is incident to a $(Q_p^d,T)$-surplus edge. In the latter case, both endpoints of that edge and the tree path between them remain.

Finally, each component spanned by $E(T)\setminus E(T^-)$ is a bare path in $C_2$ and has as an endpoint a distinct leaf of $T^-$. Since each such path has at most $d^3-3$ vertices we have that
\[
 |V(T)|\leq |V(T^-)|+(d^3-3)\ell(T^-)\leq d^3|V(T^-)|.
\]

To obtain $T^+$, start from $T$ and proceed in two stages. First, whenever a leaf has a degree-two parent, is not incident to a surplus edge, and has degree two in $C_2$, append its other edge in $C_2$. Repeat until no such leaf remains. This stage preserves the number of leaves.

Next, consider a leaf that is not incident to a surplus edge and whose parent does not have degree at least three. After the first stage, such a leaf has degree at least three in $C_2$. Append two of its edges in $C_2$ that are absent from the current tree. Their other endpoints lie outside the tree, since the leaf is not incident to a surplus edge. Repeat until no such leaf remains, and denote the resulting tree by $T^+$. This stage increases the number of leaves by at most a factor of two, so $\ell(T)\leq\ell(T^+)\leq 2\ell(T)$.%

By construction, $T^+$ is special. Every non-leaf vertex of $T^+$ that is neither a parent of a leaf nor a vertex of $T$ has degree two in $C_2$. Since $T$ is interior-closed, it follows that $T^+$ is almost interior-closed. 

Finally, each component spanned by $E(T^+)\setminus E(T)$ is a bare path in $C_2$, possibly extended by two edges, and has at most $d^3-1$ vertices. Consequently,
\[
 |V(T^+)|\leq |V(T)|+(d^3-1)\ell(T)\leq d^3|V(T)|,
\]
as required.
\end{proof}

\subsubsection{Special trees in $Q_p^d$}
For the rest of this section let $\cE$ be the event that every connected set $U$ in $Q_p^d$ with $ m/d^4\le |U|\le Md^6$ spans at most $(1+\eps^{12})|U|$ edges in $Q_p^d$.
\begin{lemma}\label{lem: event E}
    Whp the event $\cE$ occurs.
\end{lemma}
\begin{proof} 
If $\cE$ does not occur, then there is a tree $T\subseteq Q_p^d$ of order $s\in[m/d^4,Md^6]$ that spans $\eps^{12}s$ additional edges in $Q_p^d$. Harper's inequality implies that a set of size $s\leq d^6M$ spans at most  $ s \log_2 s\leq 2\alpha ds$ edges in $Q^d$. Thus Lemma~\ref{lem: tree count} gives
\begin{align*}
\Pr(\neg \cE)&\leq  \sum_{s=m/d^4}^{d^6M}2^d(ed)^{s-1}p^{s-1}  \binom{2\alpha  ds}{{\eps^{12}s}}p^{{\eps^{12}s} } \leq \sum_{s=m/d^4}^{d^6M}2^d3^s\bfrac{2 e\alpha  d p s}{{\eps^{12}s}}^{{\eps^{12}s}} 
\\&\leq \sum_{s=m/d^4}^{d^6M}2^d3^s\alpha^{{\eps^{13}s}}= \sum_{s=m/d^4}^{d^6M}2^d3^s\eps ^{{s}/{\eps^{7}}}=o(1). 
\end{align*}
\end{proof}

We now bound the numbers of leaves of special trees in $C_2$. For a tree of order $s$ in the range below, the bound on $b$ shows that few leaves are supported by surplus edges. The bounds on $a$ then give $O(\eps s)$ leaves in total, and $\Omega(\eps s)$ leaves when the tree is almost interior-closed.

\begin{lemma}\label{lem: existance of special trees}
Let
\[\alpha=\alpha(\eps)\coloneqq \eps^{1/\eps^{20}}, \quad m=m(\eps,d) \coloneqq \eps^2\cdot  2^{\alpha d},  \quad \text{and} \quad M=M(\eps,d)\coloneqq \eps^{-2}\cdot  2^{\alpha d}.\]
Whp every $Q^d_p$-special $(a,b)$-tree $T\subseteq C_2$ of order $s$ with $ m/d^4\le s\le Md^5$ 
satisfies
\[
 b<\eps^{10}s,\qquad a\le1000\eps s.
\]
If, in addition, $T$ is almost interior-closed in $C_2$, then
$  a\ge\eps s/1000$.
\end{lemma}

\begin{proof}
Assume that there exists a $Q^d_p$-special $(a,b)$-tree $T\subseteq C_2$ of order $s$ with $ m/d^4\le s\le Md^5$.  The proof is divided into $4$ cases, which we rule out. Cases 1 and 2 give the upper bounds on $b$ and $a$, respectively. For Cases 3 and 4, we use almost interior-closedness to prove the lower bound on $a$.
\vspace{3mm}
\\ \noindent \textbf{Case 1: $b \geq \eps^{10}s$.} 
Each leaf in $B(T)$ meets a surplus edge, and a surplus edge meets at most two such
leaves. Hence $V(T)$ spans at least $|V(T)|-1+b/2\geq (1+\eps^{11})|V(T)|$ edges, contradicting the event $\cE$. Hence Case 1 does not occur whp.
\medskip\noindent
\\For the remaining cases, we first record a necessary condition on the leaves of $T$. For a leaf $v\in V(T)$ we denote by $C(v,T)$ the component of $v$ in $Q_p^d-(V(T)\setminus \{v\})$.
\vspace{3mm}
\\ \noindent \textbf{Property (P1)}: For every leaf vertex $v$ of $T$ at least one of the following holds. 
\begin{itemize}
    \item $C(v,T)$ has size at least $\eps^{-6}$.
    \item $C(v,T)$ has size less than $\eps^{-6}$ and it spans a cycle.
    \item $v$ is connected to another vertex of $T$ via a path of length at most $\eps^{-6}$ that uses no edge of $T$.
\end{itemize}
(P1) follows from the fact that  $V(T)\subseteq V(C_2)$ implies that every vertex $v\in V(T)$ either lies on a cycle or on a path connecting two cycles in $Q_p^d$. To bound the probability that a given tree $T$ satisfies (P1), we use the following exploration procedure.

Let $r=\eps^{-6}$ and enumerate the leaves as $L(T)=\{v_1,\ldots,v_{a+b}\}$. We explore from them in this order. Having constructed the exploration trees $F_1,\ldots,F_{i-1}$, let $U_i=V(T)\cup\bigcup_{j<i}V(F_j)$ and start $F_i$ at $v_i$. Extend $F_i$ by depth-first search, querying edges from its vertices to $V(Q^d)\setminus((U_i\setminus\{v_i\})\cup V(F_i))$. Stop when $F_i$ reaches $r$ vertices or when the component of $v_i$ in $Q_p^d-(U_i\setminus\{v_i\})$ has been fully explored.

For $i\in [a+b]$, call the leaf $v_i$ successful if $|V(F_i)|=r$. Note that conditioned on $F_1,\dots F_{i-1}$, the exploration process for determining $F_i$ is dominated by a branching process with offspring distribution  $\operatorname{Binomial}(d,p)$.  Thus, by \eqref{eq: exploring to medium}, conditioned on $F_1,\dots F_{i-1}$ the vertex $v_i$ is successful with probability at most $2\eps+O(\eps^2)\le5\eps/2$. 
 
Let $U=V(T)\cup\bigcup_{j\leq a+b}V(F_j)$. If (P1) is satisfied, then every non-successful exploration tree $F_i$ is incident to an unqueried present edge spanned by $U$ that does not belong to $T\cup\bigcup_jF_j$. Since this union is a spanning tree on $U$, the number $z$ of such surplus edges satisfies, in the event $\mathcal E$,
\[
 z\leq\eps^{11}|U|+1\leq\eps^{11}rs+1\leq 2\eps^5s.
\]
Here $|U|\leq rs\leq Md^6$. Each surplus edge is incident to at most two exploration trees, so at most $4\eps^5s$ leaves are non-successful.  If $a\geq\eps^3s$ and $b\leq\eps^{10}s$, the probability that at least $a-4\eps^5s$ leaves are successful is at most
\begin{align*}
 \binom{a+b}{b+4\eps^5s}(5\eps/2)^{a-4\eps^5s}
 \leq (5\eps/2)^a\exp\bigl(O(\eps^2\log(1/\eps))a\bigr)\leq(3\eps)^a.
\end{align*}
Let $\mathcal H_1$ be the full exploration history and let $\mathcal A_1(T)$ be the event that at least $a-4\eps^5s$ leaves are successful. This event is determined by $\mathcal H_1$, and $\mathcal E$ together with (P1) implies $\mathcal A_1(T)$. The exploration does not query edges of $T$. Therefore, for every fixed tree $T\subseteq Q^d$ in this range,
\begin{align}\label{eq: bound a-related}
\Pr(\text{the event }\cE\text{ occurs and }\text{ $T$ satisfies (P1)})\leq  \Pr\bigl( \mathcal A_1(T)\bigr)
 \leq (3\eps)^a.
\end{align}
\noindent\medskip\noindent\textbf{Case 2: $b\leq\eps^{10}s$ and $a>1000\eps s$.}
If $a+b>s/4$, then $a>s/5$ for sufficiently small $\eps$. Lemma~\ref{lem: tree count} and \eqref{eq: bound a-related} bound the probability of such a tree, together with $\mathcal E$, by
\[
 \sum_{s=m/d^4}^{d^5M}2^d(edp)^{s-1}(3\eps)^{s/5}=o(1).
\]
Otherwise $s\geq 4(a+b)$ and $a>b$, and the tree has a vertex of degree at least three at which to root it. Lemma~\ref{lem: count prune trees} and \eqref{eq: bound a-related} now give the following bound for the remaining part of Case 2:
\begin{align*}
 &\Pr(\neg\mathcal E)+
 \sum_{s=m/d^4}^{d^5M}\sum_{a=1000\eps s}^{s}\sum_{b=0}^{\eps^{10}s}
 2^d(a+b)d^{s-1}\left(\frac{2e^3s}{a}\right)^{a+2b}p^{s-1}(3\eps)^a\\
 &\quad\leq o(1)+\sum_{s=m/d^4}^{d^5M}\sum_{a=1000\eps s}^{s}
 2^{3d}(pd)^{s-1}
 \left(\frac{6\eps e^3s}{a}\right)^a
 \left(\frac{2e^3s}{a}\right)^{2\eps^{10}s}\\
 &\quad\leq o(1)+\sum_{s=m/d^4}^{d^5M}\sum_{a=1000\eps s}^{s}
 2^{3d}e^{\eps s}
 \left(\frac{6e^3}{1000}\cdot(2e^3\eps^{-1})^{\eps^9}\right)^a
 =o(1).
\end{align*}
Hence Case 2 does not occur whp.
\vspace{3mm}
\\ \noindent
For cases 3 and 4, we use the additional restriction imposed by almost interior-closedness. Let $U$ and $\mathcal H_1$ be as defined at the end of the preceding exploration. Let $K(T)$ be the set of vertices of $T$ that are neither leaves nor parents of leaves. For $v\in K(T)$ let $C(v,T)$ be the component of $v$ in $Q_p^d-(U\setminus\{v\})$.

\medskip\noindent\textbf{Property (P2):} For $v\in K(T)$, $C(v,T)$ is a tree.

Indeed, if (P2) does not hold, then $v$ has a neighbor outside $U$ that also belongs to $C_2$. This contradicts the assumption that $T$ is almost interior-closed in $C_2$. Let $\mathcal A_2(T)$ be the event that $T$ satisfies (P2). Let $\mathcal E_2$ be the high-probability event that every component of $Q_p^d-V(C_2)$ has at most $d^3$ vertices, as in Corollary~\ref{cor: sizes of paths bare paths}.

Conditional on $\mathcal H_1$, we couple a further exploration with $k=|K(T)|$ independent branching processes with offspring distribution $\operatorname{Binomial}(D,p)$, where $D= (1-\eps^4)d$. The coupling will be such that, in the event $\mathcal E_2\cap\mathcal A_2(T)$, every descendant corresponds to a distinct newly discovered vertex outside $C_2$. In particular, all $k$ processes die out.

Let $R=K(T)$, $R^+=U$, and let $F$ consist of the unqueried edges in $Q^d$ from $R$ to $V(Q^d)\setminus R^+$. Sequentially, for $i=1,\ldots,k$, start a queue $Q_i=\{u_i\}$ with one artificial root. While $Q_i\neq\varnothing$, remove a vertex from $Q_i$ and query $D$ edges from $F$, one at a time. After each query, remove that edge from $F$. If it is present, let $a_e$ be its endpoint outside $R^+$, add $a_e$ to both $R$ and $R^+$, and add one child to $Q_i$. Update $F$ immediately by adding the unqueried edges from $a_e$ to outside $R^+$ and removing edges now spanned by $R^+$. Thus each present query discovers a new vertex. If no eligible edge remains before the $D$ queries have been completed, return FAILURE and complete the branching processes with fresh independent Bernoulli trials. Terminate the $i$th iteration once $Q_i=\varnothing$.

Now assume that the event $\mathcal E_2\cap\mathcal A_2(T)$ occurs and $a+b\leq1001\eps s$. We will show that FAILURE does not occur. Indeed, every discovered vertex outside $K(T)$ lies in a tree attached to a vertex of $K(T)$, so there are at most $d^4k$ such vertices. Thus,  $|R^+| \leq |U|+d^4 |K(T)|\leq 2d^4s$. Therefore $|R|/|R^+|\geq (s-2(a+b))/rs\geq 1/2r$. Harper's inequality gives that $R^+$ spans at most $\alpha^{2/3}d|R^+|\leq \eps^{10}d|R|$ edges. Hence $e_{Q^d}(R, V(Q^d)\setminus R^+)\geq (1-\eps^5)d|R|$. If FAILURE occurs, then it occurs when the $t$-th artificial vertex is processed for some $t\leq |R|$ at which point $F=\emptyset$. Thus at most $Dt$ edges had entered $F$ and had been processed. On the other hand, each of the at least $(1-\eps^5)d|R|$ edges from $R$ to $V(Q^d)\setminus R^+$ had entered $F$ and had been processed. Since $(1-\eps^5)d|R|>D|R|$ we have that FAILURE cannot occur. 
 
In particular, all of the $k$ processes die out. 
The extinction probability of each process is $1-2\eps+O(\eps^2)$. Uniformly over the full history $\mathcal H_1$, we obtain, for $a+b\leq1001\eps s$,
\begin{align}\label{eq: P2-associated bound}
 &\Pr\brac{\mathcal E_2\cap\mathcal A_2(T)
       \mid \mathcal H_1} \leq (1-2\eps+O(\eps^2))^k\leq e^{-7\eps k/4}\leq e^{-3\eps s/2}.
\end{align}

\medskip\noindent\textbf{Case 3: $T$ is almost interior-closed in $C_2$, $b\leq\eps^{10}s$ and $a<\eps^3s$.} If  $T$ has a vertex of degree at least $3$, then it has at least $3$ leaves, hence $b\geq \max\{0,3-a\}$, and Lemma~\ref{lem: count prune trees} applies. Otherwise $T$ is a path.  Lemma~\ref{lem: count prune trees} and \eqref{eq: P2-associated bound} give the following upper bound on the probability of Case 3:
\begin{align*}
&\Pr(\neg \cE_2)+\sum_{s=m/d^4}^{d^5M}   2^d d^{s-1}  p^{s-1}      e^{-3\eps s/2}\\&+\sum_{s=m/d^4}^{d^5M} \sum_{a=0}^{\eps^3 s} \sum_{b=\max\{0,3-a\}}^{\eps^{10}s}  2^d (a+b) d^{s-1}\bfrac{2e^3s}{a+2b}^{a+2b}  p^{s-1}      e^{-3\eps s/2}
\\&\leq  \sum_{s=m/d^4}^{d^5M}   2^{4d} (dp)^{s}\bfrac{2e^3}{\eps^3 }^{3\eps^3 s}       e^{-3\eps s/2} \leq  \sum_{s=m/d^4}^{d^5M}   2^{4d} \brac{\bfrac{2e^3}{\eps^3 }^{3\eps^2 }       e^{-1/2}}^{\eps s}=o(1).   
\end{align*}

\medskip\noindent\textbf{Case 4: $T$ is almost interior-closed in $C_2$, $b\leq\eps^{10}s$ and $\eps^3s\leq a<\eps s/1000$.}
Lemma~\ref{lem: count prune trees}, \eqref{eq: bound a-related} and \eqref{eq: P2-associated bound}  give the following upper bound on the probability of Case 4:
\begin{align*}
 &\Pr(\neg\mathcal E)+\Pr(\neg\mathcal E_2) +\sum_{s=m/d^4}^{d^5M}\sum_{a=\eps^3s}^{\eps s/1000}
 \sum_{b=0}^{\eps^{10}s}2^d(a+b)d^{s-1}
 \left(\frac{2e^3s}{a}\right)^{a+2b}p^{s-1}(3\eps)^a e^{-3\eps s/2}\\
 &\leq o(1)+\sum_{s=m/d^4}^{d^5M}\sum_{a=\eps^3s}^{\eps s/1000}
 2^{3d}\left(\frac{6e^3\eps s}{a}\right)^a
 \left(\frac{2e^3s}{a}\right)^{2\eps^7a}e^{-\eps s/2}\\
 &\leq o(1)+\sum_{s=m/d^4}^{d^5M}\sum_{a=\eps^3s}^{\eps s/1000}
 2^{3d}\left(\frac{7e^3\eps s}{a}\right)^a e^{-\eps s/2}\\
 &\leq o(1)+\sum_{s=m/d^4}^{d^5M}\sum_{a=\eps^3s}^{\eps s/1000}
 2^{3d}\left((7000e^3)^{1/1000}e^{-1/2}\right)^{\eps s}
 =o(1).
\end{align*}

 \end{proof}

\begin{proof}[\textbf{Proof of Theorem~\ref{thm: expansion of 2-core}}]
Let $\alpha,m,M$ be as in Theorem~\ref{thm: expansion of 2-core}. Let $\mathcal E^+$ be the event that $\cE$  and the high-probability events in Lemmas~\ref{lem: special trees} and~\ref{lem: existance of special trees} occur. It suffices to prove that in the event $\cE^+$, (1) and (2) hold. For the lower bound in (1), we extend a tree witness and use the lower bound on its number of leaves. For the upper bound we prune a tree witness and use the upper bound on its number of leaves.
For (2), the surplus-edge bound suffices.  

First, suppose that the lower bound in (1) does not hold, and let $S\subseteq V(C_2)$ be connected and of size $|S|\in [m,M+1]$ such that $|N_{C_2}(S)|\leq \eps |S|/4000$. Then, $|S\cup N_{C_2}(S)|\leq 2(M+1)$ and therefore $S\cup N_{C_2}(S)$ spans at most $|S\cup N_{C_2}(S)|+t$ edges, where $t \coloneqq  \eps^{12}|S\cup N_{C_2}(S)|$.
Choose a spanning tree on $C_2[S]$, attach to it the vertices in $N_{C_2}(S)$ via the edges in $C_2$ and let $T$ be the resulting tree. Then $T$ is  interior-closed  in $C_2$. 

Every leaf of $T$ lying in $S$ is incident
to a surplus edge, since its degree in $C_2$ is at least two. There are at most $t+1$ surplus edges,
so $\ell(T) \leq |N_{C_2}(S)| + 2(t + 1)$. Lemma~\ref{lem: special trees} implies that $T$ can be extended into a $Q_p^d$- special $(a',b')$-tree $T^+$ that is almost interior-closed in $C_2$, has at most $d^3|V(T)|\leq2d^3(M+1)$ vertices, and has at most $2\ell(T)$ leaves.
Therefore,
\[ \ell(T^+)\leq2\ell(T)\leq2|N_{C_2}(S)|+4(t+1)
 <\frac{\eps|S|}{1000}\leq\frac{\eps|V(T^+)|}{1000}.
\]
This contradicts the lower bound on $a'$ given in Lemma~\ref{lem: existance of special trees}. Hence in the event $\mathcal E^+$, the lower bound in (1) holds.

Second, suppose that the upper bound in (1) does not hold. Let $S\subseteq V(C_2)$ be connected and of size $s=|S|\in[m,M+1]$ such that $|N_{C_2}(S,V(C_2)\setminus S)|>2000\eps s$. Choose a spanning tree on $C_2[S]$ and attach  all vertices of $N_{C_2}(S)$ as leaves, obtaining an interior-closed tree $T$ of order at most $(d+1)s$. Lemma~\ref{lem: special trees} implies that $T$ contains a $Q_p^d$-special $(a',b')$-tree $T^-$ with $|V(T^-)|\geq|V(T)|/d^3$ and $\ell(T^-)=\ell(T)$. Its order lies in the range of Lemma~\ref{lem: existance of special trees}, and
\[
 \frac{a'+b'}{|V(T^-)|}
 =\frac{\ell(T^-)}{|V(T^-)|}  \geq\frac{|N_{C_2}(S,V(C_2)\setminus S)|}{s+|N_{C_2}(S,V(C_2)\setminus S)|} 
 \geq\frac{2000\eps s}{(1+2000\eps)s} >1500\eps.
\]
This contradicts the bounds $a'\leq1000\eps|V(T^-)|$ and $b'<\eps^{10}|V(T^-)|$ in Lemma~\ref{lem: existance of special trees}. Hence in the event $\mathcal E^+$, the upper bound in (1) holds.

Finally, if (2) fails, there are a connected set   $S$ and $Z\subseteq N_{C_2}(S)$   such that $S\cup Z$ spans more than $(1+\eps^{12})|S\cup Z|$ edges, contradicting $\mathcal E$. Hence in the event $\mathcal E^+$, (2) holds.
\end{proof}

\subsection*{Acknowledgements}
This research was funded in whole or in part by the Austrian Science Fund (FWF) 10.55776/ ESP3863424 and the SNSF grant 200021-228014.  For open access purposes, the authors have applied a CC BY public copyright license to any author accepted manuscript version arising from this submission.

\paragraph{Statement on AI use.} 
The results presented in this paper were obtained without AI assistance. ChatGPT Pro was used for text editing, polishing, and proofreading.

\bibliographystyle{abbrv}
\bibliography{main}

\end{document}